\documentclass[12pt,oneside,reqno]{amsart}
\usepackage{}
\usepackage{amssymb}
\usepackage[colorlinks=true, linkcolor=blue, citecolor=blue]{hyperref}
\usepackage{bbm}
\usepackage{cases}
\usepackage{amsmath}
\usepackage{graphicx}
\usepackage{cite}
\usepackage{mathrsfs}
\usepackage{stmaryrd}
\usepackage{color}
\usepackage{soul}
\usepackage[dvipsnames]{xcolor}
\usepackage{amsfonts}
\usepackage{enumerate,amsmath,amssymb,amsthm}

\numberwithin{equation}{section}

\newcommand{\be}{\begin{eqnarray}}
\newcommand{\mE}{\end{eqnarray}}
\newcommand{\ce}{\begin{eqnarray*}}
\newcommand{\de}{\end{eqnarray*}}
\newtheorem{theorem}{Theorem}[section]
\newtheorem{lemma}[theorem]{Lemma}
\newtheorem{remark}[theorem]{Remark}
\newtheorem{definition}[theorem]{Definition}
\newtheorem{proposition}[theorem]{Proposition}
\newtheorem{example}[theorem]{Example}
\newtheorem{corollary}[theorem]{Corollary}

\def\e{{\mathrm{e}}}
\def\eps{\varepsilon}

\def\p{\partial}

\def\[{{\Big[}}
\def\]{{\Big]}}
\def\<{{\langle}}
\def\>{{\rangle}}
\def\({{\Big(}}
\def\){{\Big)}}

\def\bx{{\mathbf{x}}}

\def\dif{{\mathord{{\rm d}}}}

\def\no{\nonumber}
\def\={&\!\!=\!\!&}
\def\bt{\begin{theorem}}
\def\et{\end{theorem}}
\def\bl{\begin{lemma}}
\def\el{\end{lemma}}
\def\br{\begin{remark}}
\def\er{\end{remark}}

\def\bd{\begin{definition}}
\def\ed{\end{definition}}
\def\bp{\begin{proposition}}
\def\ep{\end{proposition}}
\def\bc{\begin{corollary}}
\def\ec{\end{corollary}}
\def\bx{\begin{example}}
\def\ex{\end{example}}

\def\cC{{\mathcal C}}

\def\cF{{\mathcal F}}

\def\cI{{\mathcal I}}

\def\cL{{\mathcal L}}

\def\cP{{\mathcal P}}

\def\cT{{\mathcal T}}

\def\cW{{\mathcal W}}

\def\mC{{\mathbb C}}

\def\mE{{\mathbb E}}
\def\mF{{\mathbb F}}

\def\mN{{\mathbb N}}

\def\mP{{\mathbb P}}

\def\mR{{\mathbb R}}

\def\sE{{\mathscr E}}
\def\sF{{\mathscr F}}

\def\sJ{{\mathscr J}}

\def\sL{{\mathscr L}}

\def\sP{{\mathscr P}}

\def\geq{\geqslant}
\def\leq{\leqslant}

\def\T{\mathord{{\rm Tr}}}

\allowdisplaybreaks

\begin{document}

\title{Poisson equation on Wasserstein space and diffusion approximations for  McKean-Vlasov equation}

\date{}

\author{Yun Li, Fuke Wu and Longjie Xie}

\address{Yun Li:
       Institute of Systems Science,
       Academy of Mathematics and Systems Science, Chinese Academy of Sciences, and School of Mathematics Sciences, University of Chinese Academy of Sciences, Beijing 100149, P.R.China\\
       Email: liyun@amss.ac.cn
}

\address{Longjie Xie:
	School of Mathematics and Statistics, Jiangsu Normal University,
	Xuzhou, Jiangsu 221000, P.R.China\\
	Email: longjiexie@jsnu.edu.cn
}

\address{Fuke Wu:
        School of Mathematics and Statistics,  Huazhong University of Science and Technology, Wuhan, Hubei 430074, P.R.China\\
        Email: wufuke@hust.edu.cn	
}

\thanks{
	This work is supported by NNSF of China (No. 12090011, 12071186,  11931004,  61873320).
}

\begin{abstract}
We consider the fully-coupled McKean-Vlasov equation with multi-time-scale potentials, and all the coefficients depend  on the distributions of both the slow component and the fast motion.
By studying the smoothness of the solution of the non-linear Poisson equation on Wasserstein space, we derive the asymptotic limit as well as the quantitative error estimate of the convergence for the slow process.
	Extra homogenized drift term containing derivative in the measure argument of the solution of the Poisson equation appears in the limit, which seems to be new and is unique for systems involving the fast distribution.

	\bigskip

	
\end{abstract}

\maketitle

\section{Introduction and main result}

For $d\geq 1$, let $\sP_2(\mR^{d})$ be the space of all square integrable probability measures over $\mR^d$ equipped with the Wasserstein metric, i.e.,
$$
\cW_2(\mu_1,\mu_2):=\inf_{\pi\in\cP(\mu_1,\mu_2)}\left(\int_{\mR^d\times\mR^d}|x-y|^2\pi(\dif x,\dif y)\right)^{\frac{1}{2}},\quad \forall \mu_1,\mu_2\in \sP_2(\mR^{d}),
$$
where $\cP(\mu_1,\mu_2)$ denotes the class of measures on $\mR^d\times\mR^d$ with marginal $\mu_1$ and $\mu_2$.
Consider the following multi-scale McKean-Vlasov stochastic differential equations (SDEs for short) in $\mR^{d_1}\times\mR^{d_2}$:
\begin{equation} \label{sde}
\left\{ \begin{aligned}
&\dif X^{\eps}_t =F(X^{\eps}_t,\cL_{X_t^\eps},Y^{\eps}_t,\cL_{Y^{\eps}_t})\dif t+\frac{1}{\eps}H(X^{\eps}_t,\cL_{X_t^\eps},Y^{\eps}_t,\cL_{Y^{\eps}_t})\dif t\\
&\qquad\qquad\qquad\qquad\qquad\quad\qquad+G(X^{\eps}_t,\cL_{X_t^\eps},Y^{\eps}_t,\cL_{Y^{\eps}_t})\dif W^1_t,\quad\,\,\,\, X^{\eps}_0=\xi,\\
&\dif Y^{\eps}_t =\frac{1}{\eps}c(X^{\eps}_t,\cL_{X_t^\eps},Y^{\eps}_t,\cL_{Y^{\eps}_t})\dif t+\frac{1}{\eps^2}b(\cL_{X_t^\eps},Y^{\eps}_t,\cL_{Y^{\eps}_t})\dif t\\
&\qquad\qquad\qquad\qquad\qquad\quad\qquad\,+\frac{1}{\eps}\sigma(\cL_{X_t^\eps},Y^{\eps}_t,\cL_{Y^{\eps}_t})\dif W_t^2,\quad\qquad\!\!   Y^{\eps}_0=\eta,
\end{aligned} \right.
\end{equation}
where $d_1, d_2\geq 1$, $\cL_\vartheta$ denotes the law of  a random
variable $\vartheta$, the coefficients $F, H: \mR^{d_1}\times\sP_2(\mR^{d_1})\times\mR^{d_2}\times\sP_2(\mR^{d_2})\to\mR^{d_1}$, $G: \mR^{d_1}\times\sP_2(\mR^{d_1})\times\mR^{d_2}\times\sP_2(\mR^{d_2})\to\mR^{d_1}\otimes\mR^{d_1}$, $c: \mR^{d_1}\times\sP_2(\mR^{d_1})\times\mR^{d_2}\times\sP_2(\mR^{d_2})\to\mR^{d_2}$, $b: \sP_2(\mR^{d_1})\times\mR^{d_2}\times\sP_2(\mR^{d_2})\to\mR^{d_2}$ and $\sigma: \sP_2(\mR^{d_1})\times\mR^{d_2}\times\sP_2(\mR^{d_2})\to\mR^{d_2}\otimes\mR^{d_2}$ are measurable functions, $W^1_t$, $W^2_t$ are $d_1$, $d_2$-dimensional independent standard Brownian motions defined on some probability space $(\Omega,\sF,\mP)$,  $\xi$, $\eta$ are $d_1$, $d_2$-dimensional random variables with finite moments, respectively, and the small parameter $0<\eps\ll 1$ represents the separation  of  time scales   between   the slow component $X_t^\eps$ (often being called the driving process)  and  the fast motion $Y_t^\eps$ (also being called the driven process).

\vspace{2mm}
The McKean-Vlasov equation, also known as the distribution dependent SDE or the mean-field  SDE, describes the limiting behaviour of an individual particle involving within a system of particles interacting through their empirical measure, as the size of the
population grows to infinity (the so-called
propagation of chaos). The pioneer work on such system was indicated by Kac \cite{K} in kinetic theory and  McKean \cite{M} in the study of non-linear partial differential equations (PDEs for short). So far, the McKean-Vlasov SDEs have been investigated in various aspects such as well-posedness, connection with non-linear Fokker-Planck equations, large deviation  and numerical approximation, etc, we refer the readers to \cite{BR,BRR,CGT,CD,CC2,CC,CGPS,CCD,CF,CM,DEGS,EGZ,FG,GP,GS,HW,K2,MV} among others.
Meanwhile, the presence of multiple scales arises naturally in many applications  ranging from climate modeling to chemical physics, and has been the central
topic of study in science and engineering (see \cite{PS}). In particular,  multiple scales can leads to hysteresis loops in the bifurcation diagram and induce phase transitions of certain McKean-Vlasov equation as studied in \cite{CGPS,DGPS,DP,GP}. Existing averaging results for slow-fast McKean-Vlasov SDEs can be found in \cite{BS2,BS,BS3,HLL2,RSX}. However, the coefficients of the systems considered in
all the previous works are not allowed to depend on the distribution of the fast motion.
Of particular relevance to us is the system of weakly interacting diffusions in a two-scale potential relying on the faster empirical measure considered in \cite{GGP}, the combined mean field and diffusive limit was investigated therein.

\vspace{2mm}
Our aim in this paper is to derive rigorously the homogenized limit of the non-linear system (\ref{sde}) as $\eps\to0$. One of the novelty of the McKean-Vlasov equation (\ref{sde}) lies in that, even the slow component $X_t^\eps$ has a fast varying term. This is known to be closely related to the homogenization of non-linear PDEs (see e.g. \cite{CC2,CCKW,DGO,HP}). In particular, when $F\equiv G\equiv0$ and $H(x,\mu,y,\nu)=y$, then system (\ref{sde}) reduces to the overdamped non-linear kinetic equation (see e.g. \cite{CC2,J} and the references therein). More importantly, all the coefficients in system (\ref{sde}) rely
on the distribution of the fast motion $Y_t^\eps$, which makes the situation more complicated. Such feature would then capture the system
explored in \cite{GGP}, and identify its asymptotic limit would  serve to give insight into the nature of phase transitions for some McKean-Vlasov systems as originally studied in \cite{D}. The key difficulties caused by these features are that, compared with the previous results  (e.g. \cite{BS2,BS,BS3,CGPS,GP,HLL2,RSX}), we need to control the fluctuations of the functional central limit type, seek a totally non-linear SDE as the frozen equation, and  handle the corresponding non-linear PDEs to serve as a corrector.

\vspace{2mm}

It turns out that the asymptotic limit for the system (\ref{sde})  will be given
in terms of the solution of an non-linear Poisson equation.  Namely, we need to consider the following Poisson equation on $\mR^{d_2}\times\sP_2(\mR^{d_2})$:
\begin{align}\label{po1}
\sL_0\Phi(x,\mu,y,\nu)=-H(x,\mu,y,\nu),
\end{align}
where $(x,\mu)\in\mR^{d_1}\times\sP_2(\mR^{d_1})$ are regarded as parameters, and for a test function $\varphi$, the operator $\sL_0$ is defined by
\begin{align}\label{l0}
\sL_0\varphi(y,\nu)&:=\sL_0(\mu,y,\nu)\varphi(y,\nu)\no\\
&:=\frac{1}{2} a(\mu,y,\nu)\p^2_y\varphi(y,\nu)+b(\mu,y,\nu)\cdot\p_y\varphi(y,\nu)\no\\
&+\int_{\mR^{d_2}}\Big[b(\mu,\tilde y,\nu)\cdot\p_{\nu}\varphi(y,\nu)(\tilde y)+\frac{1}{2} a(\mu,\tilde y,\nu)\cdot\p_{\tilde y}\big[\p_{\nu}\varphi(y,\nu)(\tilde y)\big]\Big]\nu(\dif \tilde y),
\end{align}
with $a(\mu,y,\nu):=\sigma\sigma^*(\mu,y,\nu)$. Though the first part of the operator involving the usual derivatives in the space variable is quite standard, the integral part contains the Lions derivative (see \cite[Section 6]{C} or \cite{L}) of the test function with respect to the measure argument.  Unlike previous works, equation (\ref{po1}) is totally non-linear. This is exactly due to the dependence on distribution of the fast motion in system (\ref{sde}). In fact,  the operator $\sL_0$ can be viewed as the infinitesimal generator of the frozen process $Y_t^{\mu,\eta}$ obtained from $Y_t^\eps$ in (\ref{sde}) by freezing the slow component at fixed $\mu\in\sP_2(\mR^{d_2})$, i.e., $Y_t^{\mu,\eta}$ satisfies the following McKean-Vlasov equation:
\begin{align}\label{sde1}
\dif Y_t^{\mu,\eta}=b(\mu,Y_t^{\mu,\eta},\cL_{Y^{\mu,\eta}_t})\dif t+\sigma(\mu,Y_t^{\mu,\eta},\cL_{Y^{\mu,\eta}_t})\dif W_t^2,\quad Y_0^{\mu,\eta}=\eta,
\end{align}
where $\mu\in\sP_2(\mR^{d_2})$ is regarded as a parameter. Furthermore, note that  the Poisson equation (\ref{po1}) is formulated on the whole space but not on compact sets (without boundary conditions). We shall show that   a necessary condition for (\ref{po1}) to be well-posed is to assume that $H$ satisfies the following centering condition (see (\ref{cen3}) below):
\begin{align}\label{cen}
\int_{\mR^{d_2}}H(x,\mu,y,\zeta^{\mu})\zeta^{\mu}(\dif y)=0,
\end{align}
where $\zeta^{\mu}(\dif y)$ is the  unique invariant measure  for the non-linear SDE (\ref{sde1}). We point out that we have freezed the $\nu$-measure variable of $H$ in (\ref{cen})  at the invariant measure $\zeta^\mu$.
Under some mild assumptions, we prove that equation (\ref{po1}) has a unique solution and admits the probabilistic representation (see {\bf Theorem \ref{po}} below)
$$
\Phi(x,\mu,y,\nu)=\mE\left(\int_0^\infty H(x,\mu,Y_t^{\mu,y,\nu},\cL_{Y_t^{\mu,\eta}}\dif t\right),
$$
where $Y_t^{\mu,y,\nu}$ satisfies the following  decoupled equation associated with the McKean-Vlasov SDE (\ref{sde1}):
\begin{align}\label{sde2}
Y_t^{\mu,y,\nu}=y+\int_0^tb\big(\mu,Y_s^{\mu,y,\nu},\cL_{Y_s^{\mu,\eta}}\big)\dif s+\int_0^t\sigma\big(\mu,Y_s^{\mu,y,\nu},\cL_{Y_s^{\mu,\eta}}\big)\dif  W^2_s
\end{align}
with $\cL_\eta=\nu$. We remark that the equation (\ref{sde2}) is not a McKean-Vlasov type SDE since the law appearing in the coefficients is not $\cL_{Y_t^{\mu,y,\nu}}$ but rather $\cL_{Y_t^{\mu,\eta}}$, i.e., the law of the solution to the McKean-Vlasov SDE (\ref{sde1}). In fact, equation (\ref{sde2}) can be viewed as a time-inhomogeneous It\^o type SDE.

\vspace{2mm}
The smoothness of the solution $\Phi$ to equation (\ref{po1}) with respect to both the $y$ and $\nu$ variables as well as the parameters $x$ and $\mu$ are also investigated. Especially, regularities with respect to the  measure $\mu$ are more difficult since $\mu$ appears not only in the function $H$ but also the McKean-Vlasov process $Y_t^{\mu,\eta}$ and its decoupled process $Y_t^{\mu,y,\nu}$. In addition, the verification that the averaged coefficients are smooth usually constitutes a separate problem connected with the smoothness of the invariant measure $\zeta^\mu$ with respect tot the parameter $\mu$ (see e.g. \cite{BSV}). Here, we  provide  explict formulas for the derivatives in the $\mu$-variable of the solution $\Phi$ as well as the averaged functions with respect to the invariant measure $\zeta^\mu$ (see Corollary \ref{avef}), which quantify the $\mu$-derivatives in terms of the $(y,\nu)$-derivatives.
Poisson equation with  parameters and in the whole space $\mR^d$ involving the classical linear  second order differential operator  was studied in a series of papers by Pardoux  and Veretennikov \cite{P-V,P-V2,P-V3}  (see also \cite{RX1})  and  is shown to be a powerful tool in
the theory of numerical approximation, stochastic averaging, homogenization  and other functional limit theorems in probability (see e.g. \cite{APV,HMP,MST,PP}).
Thus, our results for the non-linear Poisson  equation (\ref{po1}) are of independent interest.

\vspace{2mm}

To formulate the asymptotic limit of system (\ref{sde}), let us define the averaged drift and diffusion coefficients by
\begin{align}
\bar F(x,\mu)&:=\int_{\mR^{d_2}}F(x,\mu,y,\zeta^{\mu})\zeta^{\mu}(\dif y),\label{bF}\\
\overline{GG^*}(x,\mu)&:=\int_{\mR^{d_2}}GG^*(x,\mu,y,\zeta^{\mu})\zeta^{\mu}(\dif y).\label{sigmaG}
\end{align}
Again, note that the $\nu$-variable of the coefficients  in (\ref{bF}) and (\ref{sigmaG}) is freezed at the invariant measure $\zeta^\mu$.
Several extra homogenized drift and diffusion coefficients will appear. We denote by
\begin{align}
\overline{H\cdot\p_x\Phi}(x,\mu)&:=\int_{\mR^{d_2}}H(x,\mu,y,\zeta^{\mu})\cdot\p_x\Phi(x,\mu,y,\zeta^{\mu})\zeta^{\mu}(\dif y),\label{bH}\\
\overline{c\cdot\p_y\Phi}(x,\mu)&:=\int_{\mR^{d_2}}c(x,\mu,y,\zeta^{\mu})\cdot\p_y\Phi(x,\mu,y,\zeta^{\mu})\zeta^{\mu}(\dif y),\label{bc1}\\
\overline{H\cdot\Phi}(x,\mu)&:=\int_{\mR^{d_2}}H(x,\mu,y,\zeta^{\mu})\cdot\Phi(x,\mu,y,\zeta^{\mu})\zeta^{\mu}(\dif y),\label{sigma}
\end{align}
and
\begin{align}
\overline{c\cdot\p_\nu\Phi}(x,\mu,y)(\tilde x)&:=\int_{\mR^{d_2}}c(\tilde x,\mu,\tilde y,\zeta^{\mu})\cdot \p_\nu\Phi(x,\mu,y,\zeta^{\mu})(\tilde y)\zeta^{\mu}(\dif \tilde y),\label{bc2}\\
\overline{\overline{c\cdot\p_\nu\Phi}}(x,\mu)(\tilde x)&:=\int_{\mR^{d_2}}\overline{c\cdot\p_\nu\Phi}(x,\mu,y)(\tilde x)\zeta^{\mu}(\dif y)\no\\
&=\int_{\mR^{d_2}}\int_{\mR^{d_2}}c(\tilde x,\mu,\tilde y,\zeta^{\mu})\cdot \p_\nu\Phi(x,\mu,y,\zeta^{\mu})(\tilde y)\zeta^{\mu}(\dif \tilde y)\zeta^{\mu}(\dif y).\label{bc3}
\end{align}
Then we shall show that, as $\eps\to 0$, the slow component $X_t^\eps$ will converge weakly to $\bar X_t$ which satisfies the following McKean-Vlasov equation:
\begin{align}\label{ave}
\dif \bar X_t&=\bar F(\bar X_t,\cL_{\bar X_t})\dif t+\overline{H\cdot\p_x\Phi}(\bar X_t,\cL_{\bar X_t})\dif t\no\\
&\quad+\overline{c\cdot\p_y\Phi}(\bar X_t,\cL_{\bar X_t})\dif t +\tilde\mE\bigg(\overline{\overline{c\cdot\p_\nu\Phi}}(\bar X_t,\cL_{\bar X_t})(\tilde{\bar{X_t}})\bigg)\dif t\no\\
&\quad+\sqrt{\overline{GG^*}+2\overline{H\cdot\Phi}}(\bar X_t,\cL_{\bar X_t})\dif W^1_t,\qquad  \bar X_0=\xi,
\end{align}
where $\tilde{\bar{X_t}}$ is a copy of $\bar X_t$ defined on a
copy $(\tilde\Omega,\tilde\sF, \tilde{\mP})$ of the original probability space  $(\Omega,\sF, \mP)$, and $\tilde\mE$ is the expectation taken with respect to $\tilde\mP$. Meanwhile, we obtain the optimal rate of convergence in terms of $\eps$.

\vspace{2mm}
We remark that the expectation term involving derivative in the measure argument of the solution of the Poisson equation in (\ref{ave}) seems to be completely  new, which is  due to the effect of the fast distribution in the coefficients. Alternatively, this term can also be  expressed as
$$
\tilde\mE\bigg(\overline{\overline{c\cdot\p_\nu\Phi}}(\bar X_t,\cL_{\bar X_t})(\tilde{\bar{X_t}})\bigg)=\overline{\overline{\overline{c\cdot\p_\nu\Phi}}}(\bar X_t,\cL_{\bar X_t}),
$$
where
$$
\overline{\overline{\overline{c\cdot\p_\nu\Phi}}}(x,\mu):=\int_{\mR^{d_1}}\overline{\overline{c\cdot\p_\nu\Phi}}(x,\mu)(\tilde x)\mu(\dif \tilde x),
$$
thus depending only on $\bar X_t$ and its distribution $\cL_{\bar X_t}$.

\vspace{2mm}
To state our main result, we make the following assumption on the coefficients:

\vspace{2mm}
\noindent{\bf (H$^{\sigma,b}$):} there exist constants $c_2>c_1\geq 0$ such that for every $\mu\in\sP_2(\mR^{d_1})$, $y_1, y_2\in\mR^{d_2}$ and $\nu_1, \nu_2\in \sP_2(\mR^{d_2})$,
\begin{align*}
&\|\sigma(\mu,y_1,\nu_1)-\sigma(\mu,y_2,\nu_2)\|^2+2\<b(\mu,y_1,\nu_1)-b(\mu,y_2,\nu_2),y_1-y_2\>\\
&\leq c_1\cW_2(\nu_1,\nu_2)^2-c_2|y_1-y_2|^2.
\end{align*}

For brevity, the spaces of functions used in this paper will be introduced in the Notations part at the end of this section. The following is the main result of this paper.

\bt\label{main}
Let $T>0$, {\bf (H$^{\sigma,b}$)} hold and $H$ satisfy the centering condition (\ref{cen}). Assume that   $F,G,H,c\in\big(\textbf{C}_b^{4,(2,2),4,(2,2)}\cap \cC_b^{4,6,(3,3)}\big)(\mR^{d_1}\times\sP_2(\mR^{d_1})\times\mR^{d_2}\times\sP_2(\mR^{d_2}))$ and $b,\sigma\in \big(\textbf{C}_b^{(2,2),4,(2,2)}\cap \cC_b^{4,6,(3,3)}\big)(\sP_2(\mR^{d_1})\times\mR^{d_2}\times\sP_2(\mR^{d_2}))$. Then for every $\varphi\in C_b^{(3,1)}(\sP_2(\mR^{d_1}))$ and $t\in[0,T]$, we have
\begin{align*}
\big|\varphi(\cL_{X_t^\eps})-\varphi(\cL_{\bar X_t})\big|\leq C_T\,\eps,
\end{align*}
where $X_t^\eps$ and $\bar X_t$ satisfy the McKean-Vlasov equation (\ref{sde}) and (\ref{ave}), respectively, and $C_T>0$ is a constant independent of $\eps$.
\et





Interacting diffusions moving in a two-scales potential depending on the empirical measure of the  process $X_t^\eps$ was considered in \cite{GP}, and systems with potential relying on the fast empirical measure $X_t^\eps/\eps$ was studied in \cite{GGP}. As correspondence, we provide the following example.

\bx
Consider the following McKean-Vlasov type stochastic Langevin equations with two-time-scales potentials   in $\mR^d$:
\begin{align}\label{ex2}
\dif X_t^\eps=F(X_t^\eps,\cL_{X_t^\eps})\dif t+\frac{1}{\eps}H(X_t^\eps/\eps,\cL_{X_t^\eps/\eps})\dif t+\sqrt{2\beta^{-1}}\dif W_t,
\end{align}
where $\beta^{-1}:=\kappa_BT$ is the temperature ($\kappa_B$ denotes the Boltzman's constant and $T$ is the absolute temperature). Different choice of the potential functions could yield numerous models.
Letting
$$
Y_t^\eps:=X_t^\eps/\eps.
$$
Then we can rewrite the system (\ref{ex2}) as
\begin{equation} \label{ex}
\left\{ \begin{aligned}
&\dif X^{\eps}_t =F(X^{\eps}_t,\cL_{X^{\eps}_t})\dif t+\frac{1}{\eps}H(Y^{\eps}_t,\cL_{Y^{\eps}_t})\dif t +\sqrt{2\beta^{-1}}\dif W_t,\\
&\dif Y^{\eps}_t =\frac{1}{\eps}F(X^{\eps}_t,\cL_{X^{\eps}_t})\dif t+\frac{1}{\eps^2}H(Y^{\eps}_t,\cL_{Y^{\eps}_t})\dif t+\frac{\sqrt{2\beta^{-1}}}{\eps}\dif W_t.
\end{aligned} \right.
\end{equation}
System (\ref{ex}) is a particular case of (\ref{sde}). Thus, according to Theorem \ref{main}, the homogenized limit for (\ref{ex2}) shall be given by
\begin{align}\label{bex}
\dif \bar X_t=(1+c_1)\, F(\bar X_t,\cL_{\bar X_t})\dif t+ c_2\, \mE F(\bar X_t,\cL_{\bar X_t})\dif t+\sqrt{2\beta^{-1}+2 c_3}\,\dif W_t,
\end{align}
where
\begin{align*}
&c_1:=\int_{\mR^d}\p_y\Phi(y,\zeta)\zeta(\dif y),\quad c_2:=\int_{\mR^{d}}\int_{\mR^{d}} \p_\nu\Phi(y,\zeta)(\tilde y)\zeta(\dif y)\zeta(\dif \tilde y),\\
&c_3:=\int_{\mR^d}H(y,\zeta)\cdot\Phi(y,\zeta)\zeta(\dif y),
\end{align*}
and $\Phi(y,\nu)$ is the unique solution to the Poisson equation
$$
\sL_0(y,\nu)\Phi(y,\nu)=-H(y,\nu).
$$
Studying the qualitative properties of the original system (\ref{ex2}) and its homogenized limit (\ref{bex}) such as the bifurcation diagram and phase transitions would be a interesting problem, we shall carry it out in a future work.
\ex

The rest of this paper is organized as follows. Section 2 is devoted to study the non-linear Poisson equation on the whole space, regularities of the solution are obtained. In Section 3 we prepare two fluctuations lemmas, and then we give the proof of Theorem \ref{main} in Section 4. Finally, an Appendix containing the proofs of some auxiliary lemmas is provided at the end of the paper.

\vspace{5mm}
\noindent{\bf Notations.} Throughout this paper, we use the following notations.

\vspace{2mm}
 Given a function $f$ on $\mR^{d_1}\times\sP_2(\mR^{d_1})\times\mR^{d_2}\times\sP_2(\mR^{d_2})$, we define
\begin{align}\label{op1}
\sL_1f(x,\mu,y,\nu)&:=\sL_1(x,\mu,y,\nu)f(x,\mu,y,\nu)\no\\
&:=F(x,\mu,y,\nu)\cdot\p_xf(x,\mu,y,\nu)\no\\
&\quad+\frac{1}{2}\mathrm{Tr}\big(GG^*(x,\mu,y,\nu)\cdot\p^2_xf(x,\mu,y,\nu)\big),\\
\sL_2f(x,\mu,y,\nu)&:=\sL_2(x,\mu,y,\nu)f(x,\mu,y,\nu):=H(x,\mu,y,\nu)\cdot\p_xf(x,\mu,y,\nu),\label{op2}
\end{align}
and
\begin{align}
\sL_3f(x,\mu,y,\nu)&:=\sL_3(x,\mu,y,\nu)f(x,\mu,y,\nu):=c(x,\mu,y,\nu)\cdot\p_yf(x,\mu,y,\nu).\label{op3}
\end{align}
We shall write
$$
\<\phi,\nu\>:=\int_{\mR^{d_2}}\phi(y)\nu(\dif y).
$$

Let us briefly recall the derivatives with respect to the measure variable introduced by Lions. The idea is to  consider  the canonical lift of a real-valued function $f: \sP_2(\mR^d)\to\mR$ into a function $\mF: L^2(\Omega)\ni X\mapsto f(\cL_X)\in\mR$. Using the Hilbert structure of the space $L^2(\Omega)$, the function $f$ is  said to be differentiable at $\mu\in \sP_2(\mR^d)$ if its canonical lift $\mF$ is
Fr\'echet differentiable at some point $X$ with $\cL_X=\mu$. By Riezs' representation theorem, the Fr\'echet derivative $D\mF(X)$, viewed as an element of $L^2(\Omega)$, can be given by a function $\p_\mu f(\mu)(\cdot): \mR^d\mapsto\mR^d$ such that
$$
D\mF(X)=\p_\mu f(\cL_X)(X).
$$
The function $\p_\mu f(\mu)(x)$ is then called the Lions derivative ($L$-derivative for short) of $f$ at $\mu$. Similarly, we can define the higher order derivatives of $f$ at $\mu$.

\vspace{2mm}

Let $d,d_1,d_2\geq 1$, $T>0$ and $k,\ell,m\in\mN=\{0,1,2,\cdots\}$. We introduce  the following spaces of functions.

\begin{itemize}

\item The space $C_b^k(\mR^{d})$. A function $f(y)$ is in $C_b^k(\mR^{d})$ if $f$ is $k$-th continuously differentiable and all its derivatives are bounded.

  \item The space $C_b^{(k,\ell)}(\sP_2(\mR^{d}))$. A function $f(\nu)$ is in $C_b^{(k,\ell)}(\sP_2(\mR^{d}))$  if the mapping $\nu\mapsto f(\nu)$ is  $k$-times continuously $L$-differentiable, and we can find a version of $\p^k_\nu f(\nu)(\tilde y_1,\cdots,\tilde y_k)$ such that for every $\nu\in\sP_2(\mR^{d})$,  the mapping $(\tilde y_1,\cdots,\tilde y_k)\mapsto \p^k_\nu f(\nu)(\tilde y_1,\cdots,\tilde y_k)$  are in $C_b^\ell(\mR^{d}\times\cdots\times\mR^d)$.

  \item The space $C_b^{2k,(k,k)}(\mR^{d}\times\sP_2(\mR^{d}))$. A function $f(y,\nu)$ is in $C_b^{2k,(k,k)}(\mR^{d}\times\sP_2(\mR^{d}))$ if for any $\nu\in \sP_2(\mR^{d})$, the mapping $y\mapsto f(y,\nu)$ is in $C_b^{2k}(\mR^{d})$, and for any $y\in \mR^{d}$, the mapping $\nu\mapsto f(y,\nu)$ is in $C_b^{(k,k)}(\sP_2(\mR^{d}))$, and for every $1\leq i\leq k$, we can find a version of $\p^i_\nu f(y,\nu)(\tilde y_1,\cdots,\tilde y_i)$ such that the mapping $(y,\tilde y_1,\cdots,\tilde y_i)\mapsto \p^i_\nu f(y,\nu)(\tilde y_1,\cdots,\tilde y_i)$ is in $C_b^{2k-i}(\mR^d\times\cdots\times\mR^d)$.


        \item The space $C_b^{(k,\ell),2k,(k,k)}(\sP_2(\mR^{d_1})\times\mR^{d_2}\times\sP_2(\mR^{d_2}))$. A function $f(\mu,y,\nu)$ is in $C_b^{(k,\ell),2k,(k,k)}(\sP_2(\mR^{d_1})\times\mR^{d_2}\times\sP_2(\mR^{d_2}))$ if for any fixed $(y,\nu)\in \mR^{d_2}\times\sP_2(\mR^{d_2})$, the mapping $\mu\mapsto f(\mu,y,\nu)$ is in $C_b^{(k,\ell)}(\sP_2(\mR^{d_1}))$, and  for every $\mu\in\sP_2(\mR^{d_1})$,   the mapping $(y,\nu)\mapsto f(\mu,y,\nu)$ is in $C_b^{2k,(k,k)}(\mR^{d_2}\times\sP_2(\mR^{d_2}))$.

\item The space $C_b^{m,(k,\ell),2k,(k,k)}(\mR^{d_1}\times\sP_2(\mR^{d_1})\times\mR^{d_2}\times\sP_2(\mR^{d_2}))$. A function $f(x,\mu,y,\nu)$ is in $C_b^{m,(k,\ell),2k,(k,k)}(\mR^{d_1}\times\sP_2(\mR^{d_1})\times\mR^{d_2}\times\sP_2(\mR^{d_2}))$ if for any $(\mu,y,\nu)$, the mapping  $x\mapsto f(x,\mu,y,\nu)$ is in $C_b^m(\mR^{d_1})$, and for fixed $x\in \mR^{d_1}$, the mapping $(\mu,y,\nu)\mapsto f(x,\mu,y,\nu)$ is in $C_b^{(k,\ell),2k,(k,k)}(\sP_2(\mR^{d_1})\times\mR^{d_2}\times\sP_2(\mR^{d_2}))$.

  \item The space $\mC_b^{(k,\ell),2k,(k,k)}(\sP_2(\mR^{d_1})\times\mR^{d_2}\times\sP_2(\mR^{d_2}))$. A function $f(\mu,y,\nu)$ is in $\mC_b^{(k,\ell),2k,(k,k)}(\sP_2(\mR^{d_1})\times\mR^{d_2}\times\sP_2(\mR^{d_2}))$ if $f\in C_b^{(k,\ell),2k,(k,k)}(\sP_2(\mR^{d_1})\times\mR^{d_2}\times\sP_2(\mR^{d_2}))$, and we can find a version of $\p^k_\mu f(\mu,y,\nu)(\tilde x_1,\cdots,\tilde x_k)$ such that the mapping $(y,\nu)\mapsto\p^\ell_{(\tilde x_1,\cdots,\tilde x_k)}\p^k_\mu f(\mu,y,\nu)(\tilde x_1,\cdots,\tilde x_k)$ is in $C_b^{2k,(k,k)}(\mR^{d_2}\times\sP_2(\mR^{d_2}))$.

       \item The space $\mC_b^{m,(k,\ell),2k,(k,k)}(\mR^{d_1}\times\sP_2(\mR^{d_1})\times\mR^{d_2}\times\sP_2(\mR^{d_2}))$. A function $f(x,\mu,y,\nu)$ is in $\mC_b^{m,(k,\ell),2k,(k,k)}(\mR^{d_1}\times\sP_2(\mR^{d_1})\times\mR^{d_2}\times\sP_2(\mR^{d_2}))$ if $f\in C_b^{m,(k,\ell),2k,(k,k)}(\mR^{d_1}\times\sP_2(\mR^{d_1})\times\mR^{d_2}\times\sP_2(\mR^{d_2}))$, and for every $x\in\mR^{d_1}$, the mapping $(\mu,y,\nu)\mapsto f(x,\mu,y,\nu)$ is in  $\mC_b^{(k,\ell),2k,(k,k)}
     (\sP_2(\mR^{d_1})\times\mR^{d_2}\times\sP_2(\mR^{d_2}))$, the mapping $(\mu,y,\nu)\mapsto\p_x^m f(x,\mu,y,\nu)$ is in $\mC_b^{(k,\ell),2k,(k,k)}(\sP_2(\mR^{d_1})\times\mR^{d_2}\times\sP_2(\mR^{d_2}))$.

\item The space $C_b^{1,m,(k,\ell),2k,(k,k)}(\mR_+\times\mR^{d_1}\times\sP_2(\mR^{d_1})\times\mR^{d_2}\times\sP_2(\mR^{d_2}))$. A function $f(t,x,\mu,y,\nu)$ is in $C_b^{1,m,(k,\ell),2k,(k,k)}(\mR_+\times\mR^{d_1}\times\sP_2(\mR^{d_1})\times\mR^{d_2}\times\sP_2(\mR^{d_2}))$ if for any $(x,\mu,y,\nu)$, the mapping  $t\mapsto f(t,x,\mu,y,\nu)$ is in $C_b^1(\mR_+)$, and for fixed $t\in \mR_+$, the mapping $(x,\mu,y,\nu)\mapsto f(t,x,\mu,y,\nu)$ is in $C_b^{m,(k,\ell),2k,(k,k)}(\mR^{d_1}\times\sP_2(\mR^{d_1})\times\mR^{d_2}\times\sP_2(\mR^{d_2}))$. Similarly, we define the space   $\mC_b^{1,m,(k,\ell),2k,(k,k)}(\mR_+\times\mR^{d_1}\times\sP_2(\mR^{d_1})\times\mR^{d_2}\times\sP_2(\mR^{d_2}))$.

\end{itemize}

We shall mainly use the above spaces with $k=1,2,3$ and $\ell=1,2$. For simplify, we denote
\begin{align*}
&\textbf{C}_b^{4,(2,2),4,(2,2)}(\mR^{d_1}\times\sP_2(\mR^{d_1})\times\mR^{d_2}\times\sP_2(\mR^{d_2}))\\
&:=\Big(\mC_b^{4,(2,2),2,(1,1)}\cap\mC_b^{4,(1,1),4,(2,2)}\Big)(\mR^{d_1}\times\sP_2(\mR^{d_1})\times\mR^{d_2}\times\sP_2(\mR^{d_2})),
\end{align*}
and
\begin{align*}
&\cC_b^{4,6,(3,3)}(\mR^{d_1}\times\sP_2(\mR^{d_1})\times\mR^{d_2}\times\sP_2(\mR^{d_2}))\\
&:=\Big(C_b^{4,(1,3),6,(3,3)}\cap C_b^{4,(3,1),6,(3,3)}\Big)(\mR^{d_1}\times\sP_2(\mR^{d_1})\times\mR^{d_2}\times\sP_2(\mR^{d_2})).
\end{align*}
When a function $f(\mu,y,\nu)$ does not dependent on the $x$-varibale, we just view it as $f(x,\mu,y,\nu)\equiv f(\mu,y,\nu)$.

\section{Poisson equations on  Wasserstein space}

\subsection{McKean-Vlasov equations and associated elliptic equations}
Let us first consider the following McKean-Vlasov equation ({\it without parameters}) in $\mR^{d_2}$:
\begin{align}\label{eqy}
\dif Y^\eta_t=b(Y^\eta_t,\cL_{Y^\eta_t})\dif t+\sigma(Y^\eta_t,\cL_{Y^\eta_t})\dif W_t,\quad Y_0^\eta=\eta\in L^2(\Omega).
\end{align}
Without the abuse of notations, we use the same symbol $b$ and $\sigma$ to denote the drift and diffusion coefficients as in the parameterized equation (\ref{sde1}), and  for every $\varphi\in C_b^{2,(1,1)}(\mR^{d_2}\times \sP_2(\mR^{d_2}))$, we still use $\sL_0$  to denote the operator
\begin{align}\label{l00}
\sL_0\varphi(y,\nu)&:=\sL_0(y,\nu)\varphi(y,\nu)\no\\
&:=\frac{1}{2} a(y,\nu)\p^2_y\varphi(y,\nu)+b(y,\nu)\cdot\p_y\varphi(y,\nu)\no\\
&\quad+\int_{\mR^{d_2}}\Big[b(\tilde y,\nu)\cdot\p_{\nu}\varphi(y,\nu)(\tilde y)+\frac{1}{2} a(\tilde y,\nu)\cdot\p_{\tilde y}\big[\p_{\nu}\varphi(y,\nu)(\tilde y)\big]\Big]\nu(\dif \tilde y),
\end{align}
where $a(y,\nu):=\sigma\sigma^*(y,\nu)$. In this situation, the assumption {\bf ($H^{\sigma,b}$)} on the coefficients becomes

 \vspace{2mm}
\noindent{\bf ($\hat H^{\sigma,b}$):} there exist constants $c_2>c_1\geq 0$ such that for every $y_1, y_2\in\mR^{d_2}$ and $\nu_1, \nu_2\in \sP_2(\mR^{d_2})$,
\begin{align*}
\|\sigma(y_1,\nu_1)-\sigma(y_2,\nu_2)\|^2+2\<b(y_1,\nu_1)-b(y_2,\nu_2),y_1-y_2\>\leq c_1\cW_2(\nu_1,\nu_2)^2-c_2|y_1-y_2|^2.
\end{align*}

\vspace{2mm}
Let $Y_t^\eta$ be the unique strong solution of the equation (\ref{eqy}). It turns out  that the law of $Y_t^{\eta}$ only depends on $\eta$ through its distribution $\cL_\eta=\nu$. Thus, given a measure $\nu\in\sP_2(\mR^{d_2})$, it makes sense to consider $\cL_{Y_t^{\eta}}$ as a function of $\nu$ without specifying the choice of the lifted random variable $\eta$, and we can define a (non-linear) semigroup $\{P_t^*\}_{t\geq 0}$ on $\sP_2(\mR^{d_2})$ by letting
$$
P_{t}^*\nu:=\cL_{Y_{t}^{\eta}}\quad {\text {with}}\quad \cL_\eta=\nu.
$$
We say that a probability measure $\zeta$ is  an invariant measure of the McKean-Vlasov equation (\ref{eqy}) or the process $Y_t^\eta$ if
$$
P_t^*\zeta=\zeta,\quad \forall t\geq 0.
$$
Under the assumption {\bf ($\hat H^{\sigma,b}$)}, it is known that (see e.g. \cite[Theorem 3.1]{W1})  there exists a unique invariant measure $\zeta$ for the equation (\ref{eqy}). Moreover, there exist constants $C_0, \lambda_0>0$ such that for every $\nu\in\sP_2(\mR^{d_2})$,
\begin{align}\label{exp}
\|P_t^*\nu-\zeta\|_{TV}+\cW_2(P_t^*\nu,\zeta)\leq C_0\,\e^{-\lambda_0t}\,\cW_2(\nu,\zeta).
\end{align}

We then introduce, for fixed $y\in\mR^{d_2}$,  the following decoupled equation of the McKean-Vlasov equation (\ref{eqy}):
\begin{align}\label{eqyy}
\dif Y_t^{y,\nu}=b(Y_t^{y,\nu},\cL_{Y_t^{\eta}})\dif t+\sigma(Y_t^{y,\nu},\cL_{Y_t^{\eta}})\dif W_t,\quad Y_0^{y,\nu}=y\in\mR^{d_2},
\end{align}
where $\cL_\eta=\nu$.
Given a function $f(y,\nu)$ on $\mR^{d_2}\times\sP_2(\mR^{d_2})$, let us define
\begin{align}\label{Ttf}
\cT_tf(y,\nu):=\mE f(Y_t^{y,\nu},\cL_{Y_t^{\eta}}),\quad \forall t>0, y\in \mR^{d_2}, \nu\in\sP_2(\mR^{d_2}).
\end{align}
We have the following result.

\bt\label{nonpde}
Assume that $\sigma, b, f\in C_b^{2,(1,1)}(\mR^{d_2}\times\sP_2(\mR^{d_2}))$. Then for every $T>0$, the function $\cT_tf(y,\nu)$ is the unique solution in $C_b^{1,2,(1,1)}([0,T]\times\mR^{d_2}\times\sP_2(\mR^{d_2}))$ of the non-linear PDE
\begin{equation} \label{kez}
\left\{ \begin{aligned}
&\partial_t U(t,y,\nu)-\sL_0(y,\nu)U(t,y,\nu)=0,\quad \forall(t,y,\nu)\in(0,T]\times\mR^{d_2}\times\sP_2(\mR^{d_2}),\\
& U(0,y,\nu)=f(y,\nu),\\
\end{aligned} \right.
\end{equation}
where $\sL_0(y,\nu)$ is defined by (\ref{l00}).
If we further assume that {\bf ($\hat H^{\sigma,b}$)} holds, then we have
that for every $y\in\mR^{d_2}$ and $\nu\in\sP_2(\mR^{d_2})$,
\begin{align}\label{expy}
\cW_2(\cL_{Y_t^{y,\nu}},\zeta)\leq C_0(1+|y|+\cW_2(\nu,\delta_0))\e^{-\lambda_0 t},
\end{align}
where $C_0, \lambda_0>0$ are positive constants independent of $t$.

In particular, $\zeta$ is also the unique invariant for the decoupled equation (\ref{eqyy}), and we have
\begin{align}\label{expty}
|\cT_tf(y,\nu)-\<f(\cdot,\zeta),\zeta\>|\leq C_0(1+|y|+\cW_2(\nu,\delta_0))\e^{-\lambda_0 t}.
\end{align}

\et

\begin{proof}
The assertion that $\cT_tf(y,\nu)$ defined by (\ref{Ttf}) is in $C_b^{1,2,(1,1)}([0,T]\times\mR^{d_2}\times\sP_2(\mR^{d_2}))$ and is the unique solution of the equation (\ref{kez}) follows by \cite[Theorem 7.2]{BLPR}.
To prove (\ref{expy}), by It\^o's formula we have
\begin{align*}
\dif \mE|Y_t^\eta-Y_t^{y,\nu}|^2&=\mE\big[2\<Y_t^\eta-Y_t^{y,\nu}, b(Y_t^{\eta},\cL_{Y^{\eta}_t})-b(Y_t^{y,\nu},\cL_{Y^{\eta}_t})\>\\
&\quad+\|\sigma(Y_t^{\eta},\cL_{Y^{\eta}_t})-\sigma(Y_t^{y,\nu},\cL_{Y^{\eta}_t})\|^2\big]\dif t\\
&\leq -c_2\mE|Y_t^\eta-Y_t^{y,\nu}|^2\dif t.
\end{align*}
By comparsion theorem we get
$$
\cW_2(\cL_{Y_t^{y,\nu}},\cL_{Y_t^\eta})\leq C_0(1+|y|+\cW_2(\nu,\delta_0))\e^{-\lambda_0 t}.
$$
This together with estimate (\ref{exp}) implies (\ref{expy}).
Moreover, we deduce that
\begin{align*}
\cT_tf(y,\nu)-\<f(\cdot,\zeta),\zeta\>&=\mE f(Y_t^{y,\nu},\cL_{Y_t^{\eta}})-\<f(\cdot,\zeta),\zeta\>\\
&=\mE f(Y_t^{y,\nu},\cL_{Y_t^{\eta}})-\mE f(Y_t^{y,\nu},\zeta)+\mE f(Y_t^{y,\nu},\zeta)-\<f(\cdot,\zeta),\zeta\>\\
&=:I_1+I_2.
\end{align*}
For the first term, we have by (\ref{exp})  that
$$
|I_1|\leq \cW_2(\cL_{Y_t^{\eta}},\zeta)\leq C_0\,\e^{-\lambda_0 t}\cW_2(\nu,\zeta).
$$
For the second term, we have by (\ref{expy}) that
\begin{align*}
|I_2|\leq C_0\cW_2(\cL_{Y_t^{y,\nu}},\zeta)\leq C_0\e^{-\lambda_0 t}(1+|y|+\cW_2(\nu,\delta_0)).
\end{align*}
The proof is finished.
\end{proof}

Fix $\lambda>0$ below. Consider the following non-linear elliptic equation:
\begin{align}\label{ellip}
\lambda U_\lambda(y,\nu)-\sL_0(y,\nu)U_\lambda(y,\nu)=f(y,\nu),
\end{align}
where $\sL_0(y,\nu)$ is defined by (\ref{l00}).
We have:

\bt\label{22}
Assume that $\sigma, b, f\in C_b^{2,(1,1)}(\mR^{d_2}\times\sP_2(\mR^{d_2}))$. Then  the function
\begin{align}\label{u1}
U_\lambda(y,\nu):=\int_0^\infty\e^{-\lambda t}\cT_tf(y,\nu)\dif t
\end{align}
is the unique solution in $C_b^{2,(1,1)}(\mR^{d_2}\times\sP_2(\mR^{d_2}))$ to the equation (\ref{ellip}).
\et
\begin{proof}
The assertion  that $U_\lambda(y,\nu)$ defined by $(\ref{u1})$ is in $C_b^{2,(1,1)}(\mR^{d_2}\times\sP_2(\mR^{d_2}))$ and solves equation (\ref{ellip}) follows by Theorem \ref{nonpde} and the integral by part formula. In fact, we have
\begin{align*}
\sL_0U_\lambda(y,\nu)&=\int_0^\infty\e^{-\lambda t}\sL_0\cT_tf(y,\nu)\dif t\\
&=\int_0^\infty\e^{-\lambda t}\p_t\cT_tf(y,\nu)\dif t=\e^{-\lambda t}\cT_tf(y,\nu)\Big|_{0}^{\infty}+\lambda\int_0^\infty\e^{-\lambda t}\cT_tf(y,\nu)\dif t\\
&=-f(y,\nu)+\lambda U_\lambda(y,\nu),
\end{align*}
which implies the desired result.
To prove the uniqueness, let $U_\lambda(y,\nu)\in C_b^{2,(1,1)}(\mR^{d_2}\times\sP_2(\mR^{d_2}))$ be a solution to equation (\ref{ellip}). Then  by It\^o's formula (see e.g. \cite[Proposition 2.1]{CF}) we have
\begin{align*}
U_\lambda(Y_t^{y,\nu},\cL_{Y_t^\eta})=U_\lambda(y,\nu)+\int_0^t\sL_0U_\lambda(Y_s^{y,\nu},\cL_{Y_s^\eta})\dif s+M_t,
\end{align*}
where  $M_t$ is a martingale given by
$$
M_t:=\int_0^t(\p_yU_\lambda)(Y_s^{y,\nu},\cL_{Y_s^\eta})\cdot\sigma(Y_s^{y,\nu},\cL_{Y_s^\eta})\dif W_s.
$$
By the product formula, we further have
\begin{align*}
\e^{-\lambda t}U_\lambda(Y_t^{y,\nu},\cL_{Y_t^\eta})&=U_\lambda(y,\nu)+\int_0^t\e^{-\lambda s}\sL_0U_\lambda(Y_s^{y,\nu},\cL_{Y_s^\eta})\dif s+\int_0^t\e^{-\lambda s}\dif M_s\\
&\quad-\lambda\int_0^t\e^{-\lambda s}U_\lambda(Y_s^{y,\nu},\cL_{Y_s^\eta})\dif s\\
&=U_\lambda(y,\nu)-\int_0^t\e^{-\lambda s}f(Y_s^{y,\nu},\cL_{Y_s^\eta})\dif s+\int_0^t\e^{-\lambda s}\dif M_s,
\end{align*}
where in the last equality we have used (\ref{ellip}). Taking expectation from both sides and letting $t\to\infty$, we obtain
$$
U_\lambda(y,\nu)=\mE\left(\int_0^\infty\e^{-\lambda s}f(Y_s^{y,\nu},\cL_{Y_s^{\eta}})\dif s\right),
$$
which implies the desired result.
\end{proof}

Now, we consider the following non-linear Poisson equation in the whole space $\mR^{d_2}\times\sP_2(\mR^{d_2})$:
\begin{align}\label{ellip2}
\sL_0(y,\nu)U(y,\nu)=-f(y,\nu).
\end{align}
This equation can be viewed as the limit $\lambda\to0$ of equation (\ref{ellip}). However, in general the integral in (\ref{u1}) is not well-defined for measurable function $f$ if $\lambda=0$. It turns out that a necessary condition for equation (\ref{ellip2}) to be well-defined is to assume that $f$ satisfies the following centering condition:
\begin{align}\label{cen3}
\int_{\mR^{d_2}}f(y,\zeta)\zeta(\dif y)=0,
\end{align}
where $\zeta$ is the unique invariant measure of the McKean-Vlasov SDE (\ref{eqy}).
In fact, if $U\in C_b^{2,(1,1)}(\mR^{d_2}\times\sP_2(\mR^{d_2}))$ satisfies the equation (\ref{ellip2}), then by It\^o's formula  we can deduce that
\begin{align*}
U(Y_t^\eta,\cL_{Y_t^\eta})=U(\eta,\cL_{\eta})+\int_0^t\sL_0 U(Y_s^\eta,\cL_{Y_s^\eta})\dif s +\int_0^t\p_yU(Y_s^\eta,\cL_{Y_s^\eta})\sigma(Y^\eta_s,\cL_{Y^\eta_s})\dif W_s,
\end{align*}
where $\sL_0$ is defined by (\ref{l00}). Letting the initial distribution $\cL_\eta$ to be the invariant measure $\zeta$ of the equation (\ref{eqy}), then we have that for every $t>0$, $\cL_{Y_t^\eta}=\zeta$. Taking expectation from both sides of the above equality, we arrive at
\begin{align*}
\int_{\mR^{d_2}}U(y,\zeta)\zeta(\dif y)=\int_{\mR^{d_2}}U(y,\zeta)\zeta(\dif y)+\int_0^t\int_{\mR^{d_2}}\sL_0 U(y, \zeta)\zeta(\dif y)\dif s,
\end{align*}
which together with equation (\ref{ellip2}) implies (\ref{cen3}).
Note that the  $\nu$-variable of the function $f(y,\nu)$ in (\ref{cen3}) is fixed at the invariant measure $\zeta$.

We have the following result.

\bt\label{pot}
Let {\bf ($\hat H^{\sigma,b}$)} hold and $f$ satisfy the centering condition (\ref{cen3}).

 \vspace{1mm}
\noindent (i) If a function $U\in C_b^{2,(1,1)}(\mR^{d_2}\times\sP_2(\mR^{d_2}))$ satisfies the equation (\ref{ellip2}), then $U$ admits the probabilistic representation
\begin{align}\label{u2}
U(y,\nu)=\int_0^\infty\cT_tf(y,\nu)\dif t,
\end{align}
where $\cT_tf$ is defined by (\ref{Ttf}).

 \vspace{1mm}
\noindent (ii) Assume that  $\sigma, b, f\in C_b^{2,(1,1)}(\mR^{d_2}\times\sP_2(\mR^{d_2}))$. Then the function $U$ defined by (\ref{u2})
is the unique solution in $C_b^{2,(1,1)}(\mR^{d_2}\times\sP_2(\mR^{d_2}))$ to the non-linear Poisson equation (\ref{ellip2}), which also satisfies the centering condition.

 \vspace{2mm}
\noindent (iii) If we further assume that $\sigma, b, f\in C_b^{2k,(k,k)}(\mR^{d_2}\times\sP_2(\mR^{d_2}))$ with $k\geq 1$, then we have $U(y,\nu)\in C_b^{2k,(k,k)}(\mR^{d_2}\times\sP_2(\mR^{d_2}))$.

\et

\begin{proof}
We divide the proof into four steps.

\vspace{1mm}
\noindent (i) Let us first show that the integral in (\ref{u2}) is well-defined. In fact, since $f$ satisfies the centering condition (\ref{cen3}), we have by (\ref{expty}) that
$$
|\cT_tf(y,\nu)|\leq C_0(1+|y|+\cW_2(\nu,\delta_0))\e^{-\lambda_0 t}.
$$
The conclusion in (i) follows by It\^o's formula and the same argument as in the proof of Theorem \ref{22}. Taking the limit $\lambda\to0$ in (\ref{u1}) we obtain
$$
\lim_{\lambda\to0}|U_\lambda(y,\nu)-U(y,\nu)|=0.
$$
Meanwhile, by Fubini's theorem and the fact that $\zeta$ is also the
unique invariant for $Y_t^{y,\nu}$ (see (\ref{expy})), we have
\begin{align*}
\int_{\mR^{d_2}}U(y,\zeta)\zeta(\dif y)&=\int_{\mR^{d_2}}\int_0^\infty\mE f(Y_t^{y,\zeta},\zeta)\dif t\zeta(\dif y)\\
&=\int_0^\infty\int_{\mR^{d_2}}\mE f(Y_t^{y,\zeta},\zeta)\zeta(\dif y)\dif t=\int_0^\infty\int_{\mR^{d_2}} f(y,\zeta)\zeta(\dif y)\dif t=0.
\end{align*}
Thus $U$ satisfies the centering condition (\ref{cen3}).

\vspace{1mm}
\noindent (ii)
Next, let $U$ be defined by (\ref{u2}), we first consider the regularity of $U$ with respect to the  $y$-variable. In view of (\ref{Ttf}) and (\ref{u2}), we deduce that
\begin{align*}
\p_y U(y,\nu)=\int_0^\infty\mE\left[\p_yf(Y^{y,\nu}_t,\cL_{Y^\eta_t})\cdot\p_y Y^{y,\nu}_t\right]\dif t,
\end{align*}
where $\p_y Y^{y,\nu}_t$ satisfies
\begin{align*}
\dif \p_y Y^{y,\nu}_t=\p_y b(Y^{y,\nu}_t,\cL_{Y^\eta_t})\cdot\p_y Y^{y,\nu}_t\dif t+\p_y \sigma(Y^{y,\nu}_t,\cL_{Y^\eta_t})\cdot\p_y Y^{y,\nu}_t \dif W_t.
\end{align*}
According to Lemma \ref{Y1} in the Appendix, we have
\begin{align*}
\sup_{y\in\mR^{d_{2}},\nu\in\sP_{2}(\mR^{d_{2}})}\mE\|\p_y Y^{y,\nu}_t\|^2\leq C_{0}\,\e^{-c_2t},
\end{align*}
which together with the boundedness of $\|\p_yf\|$ yields
\begin{align*}
\sup_{y\in\mR^{d_{2}},\nu\in\sP_{2}(\mR^{d_{2}})}\|\p_y U(y,\nu)\|\leq C_{0}<\infty.
\end{align*}
Similarly, we have
\begin{align*}
\p^2_y U(y,\nu)=\int_0^\infty\mE\left[\p^2_yf(Y^{y,\nu}_t,\cL_{Y^\eta_t})\cdot\p_y Y^{y,\nu}_t\cdot\p_y Y^{y,\nu}_t+\p_yf(Y^{y,\nu}_t,\cL_{Y^\eta_t})\cdot\p^2_y Y^{y,\nu}_t\right]\dif t,
\end{align*}
where $\p^2_y Y^{y,\nu}_t$ satisfies
\begin{align*}
\dif \p^2_y Y^{y,\nu}_t&=\p^2_y b(Y^{y,\nu}_t,\cL_{Y^\eta_t})\cdot\p_y Y^{y,\nu}_t\cdot\p_y Y^{y,\nu}_t\dif t+\p_y b(Y^{y,\nu}_t,\cL_{Y^\eta_t})\cdot\p^2_y Y^{y,\nu}_t\dif t\\
&\quad+\p^2_y \sigma(Y^{y,\nu}_t,\cL_{Y^\eta_t})\cdot\p_y Y^{y,\nu}_t\cdot\p_y Y^{y,\nu}_t \dif W_t+\p_y \sigma(Y^{y,\nu}_t,\cL_{Y^\eta_t})\cdot\p^2_y Y^{y,\nu}_t \dif W_t.
\end{align*}
Applying Lemma \ref{Y1}  in the Appendix again we have
\begin{align*}
\sup_{y\in\mR^{d_{2}},\nu\in\sP_{2}(\mR^{d_{2}})}\|\p^2_y U(y,\nu)\|\leq C_{0}<\infty.
\end{align*}
As a result, for every $\nu\in\sP_{2}(\mR^{d_{2}})$ we obtain that $U(\cdot,\nu)\in C_b^2(\mR^{d_2})$.

\vspace{1mm}
\noindent (iii)
We proceed to study the regularity of $U$ with respect to the  $\nu$-variable.
By the definition of the $L$-derivative, we have (see e.g. \cite[Lemma 6.1]{BLPR})
\begin{align}
\p_\nu U(y,\nu)(\tilde y)&=\int_0^\infty\mE\Big[\p_yf(Y^{y,\nu}_t,\cL_{Y^\eta_t})\cdot\p_\nu Y^{y,\nu}_t(\tilde y)+\tilde\mE\big[\p_\nu f(Y^{y,\nu}_t,\cL_{Y^\eta_t})(\tilde Y^{\tilde y,\nu}_t)\cdot\p_y \tilde Y^{\tilde y,\nu}_t\big] \no\\
&\qquad\quad+\tilde\mE\big[\p_\nu f(Y^{y,\nu}_t,\cL_{Y^\eta_t})(\tilde Y^{\tilde\eta}_t)\cdot \tilde Z^{\tilde \eta}_t(\tilde y)\big]\Big]\dif t,\label{nuU}
\end{align}
and thus
\begin{align*}
\p_{\tilde y}\p_\nu U(y,\nu)(\tilde y)&=\int_0^\infty\mE\Big[\p_yf(Y^{y,\nu}_t,\cL_{Y^\eta_t})\cdot\p_{\tilde y}\p_\nu Y^{y,\nu}_t(\tilde y)\\
&\qquad\qquad+\tilde\mE\big[\p_\nu f(Y^{y,\nu}_t,\cL_{Y^\eta_t})(\tilde Y^{\tilde y,\nu}_t)\cdot\p^2_y \tilde Y^{\tilde y,\nu}_t\big] \\
&\qquad\qquad+\tilde\mE\big[\p_{\tilde y}\p_\nu f(Y^{y,\nu}_t,\cL_{Y^\eta_t})(\tilde Y^{\tilde y,\nu}_t)\cdot\p_y \tilde Y^{\tilde y,\nu}_t\cdot\p_y \tilde Y^{\tilde y,\nu}_t\big] \\
&\qquad\qquad+\tilde\mE\big[\p_\nu f(Y^{y,\nu}_t,\cL_{Y^\eta_t})(\tilde Y^{\tilde \eta}_t)\cdot\p_{\tilde y} \tilde Z^{\tilde \eta}_t(\tilde y)\big]\Big]\dif t,
\end{align*}
where $Z^\eta_t(\tilde y):=\p_\nu Y^{y,\nu}_t(\tilde y)|_{y=\eta}$ with  $\p_\nu Y^{y,\nu}_t(\tilde y)$  satisfying
\begin{align*}
\dif \p_\nu Y^{y,\nu}_t(\tilde y)&=\p_y b(Y^{y,\nu}_t,\cL_{Y^\eta_t})\cdot\p_\nu Y^{y,\nu}_t(\tilde y)\dif t+\tilde\mE\big[\p_\nu b(Y^{y,\nu}_t,\cL_{Y^\eta_t})(\tilde Y^{\tilde y,\nu}_t)\cdot\p_y \tilde Y^{\tilde y,\nu}_t\big]\dif t \\
&\quad+\tilde\mE\big[\p_\nu b(Y^{y,\nu}_t,\cL_{Y^\eta_t})(\tilde Y^{\tilde \eta}_t)\cdot \tilde Z^{\tilde \eta}_t(\tilde y)\big]\dif t+\p_y \sigma(Y^{y,\nu}_t,\cL_{Y^\eta_t})\cdot\p_\nu Y^{y,\nu}_t(\tilde y) \dif W_t\\
&\quad+\tilde\mE\big[\p_\nu \sigma(Y^{y,\nu}_t,\cL_{Y^\eta_t})(\tilde Y^{\tilde y,\nu}_t)\cdot\p_y \tilde Y^{\tilde y,\nu}_t\big]\dif W_t \\
&\quad+\tilde\mE\big[\p_\nu \sigma(Y^{y,\nu}_t,\cL_{Y^\eta_t})(\tilde Y^{\tilde\eta}_t)\cdot \tilde Z^{\tilde \eta}_t(\tilde y)\big]\dif W_t,
\end{align*}
and $\p_{\tilde y} Z^\eta_t(\tilde y):=\p_{\tilde y}\p_\nu Y^{y,\nu}_t(\tilde y)|_{y=\eta}$ with $\p_{\tilde y}\p_\nu Y^{y,\nu}_t(\tilde y)$ satisfying
\begin{align*}
\dif \p_{\tilde y}\p_\nu Y^{y,\nu}_t(\tilde y)&=\p_y b(Y^{y,\nu}_t,\cL_{Y^\eta_t})\cdot\p_{\tilde y}\p_\nu Y^{y,\nu}_t(\tilde y)\dif t+\tilde\mE\big[\p_\nu b(Y^{y,\nu}_t,\cL_{Y^\eta_t})(\tilde Y^{\tilde y,\nu}_t)\cdot\p^2_y \tilde Y^{\tilde y,\nu}_t\big]\dif t \\
&\quad+\tilde\mE\big[\p_{\tilde y}\p_\nu b(Y^{y,\nu}_t,\cL_{Y^\eta_t})(\tilde Y^{\tilde y,\nu}_t)\cdot\p_y \tilde Y^{\tilde y,\nu}_t\cdot\p_y \tilde Y^{\tilde y,\nu}_t\big]\dif t \\
&\quad+\tilde\mE\big[\p_\nu b(Y^{y,\nu}_t,\cL_{Y^\eta_t})(\tilde Y^{\tilde \eta}_t)\cdot \p_{\tilde y}\tilde Z^{\tilde \eta}_t(\tilde y)\big]\dif t\\
&\quad+\p_y \sigma(Y^{y,\nu}_t,\cL_{Y^\eta_t})\cdot\p_{\tilde y}\p_\nu Y^{y,\nu}_t(\tilde y) \dif W_t\\
&\quad+\tilde\mE\big[\p_\nu \sigma(Y^{y,\nu}_t,\cL_{Y^\eta_t})(\tilde Y^{\tilde y,\nu}_t)\cdot\p^2_y \tilde Y^{\tilde y,\nu}_t\big]\dif W_t \\
&\quad+\tilde\mE\big[\p_{\tilde y}\p_\nu \sigma(Y^{y,\nu}_t,\cL_{Y^\eta_t})(\tilde Y^{\tilde y,\nu}_t)\cdot\p_y \tilde Y^{\tilde y,\nu}_t\cdot\p_y \tilde Y^{\tilde y,\nu}_t\big]\dif W_t \\
&\quad+\tilde\mE\big[\p_\nu \sigma(Y^{y,\nu}_t,\cL_{Y^\eta_t})(\tilde Y^{\tilde\eta}_t)\cdot \p_{\tilde y}\tilde Z^{\tilde \eta}_t(\tilde y)\big]\dif W_t.
\end{align*}
By Lemma \ref{Y2} in the Appendix, we have
\begin{align*}
&\sup_{y\in\mR^{d_{2}},\nu\in\sP_{2}(\mR^{d_{2}}),\tilde y\in\mR^{d_{2}}}\mE\|\p_\nu Y^{y,\nu}_t(\tilde y)\|^2\leq C_{0}\,\e^{-(c_2-c_1-\gamma)t},
\end{align*}
and
\begin{align*}
&\sup_{y\in\mR^{d_{2}},\nu\in\sP_{2}(\mR^{d_{2}}),\tilde y\in\mR^{d_{2}}}\mE\|\p_{\tilde y}\p_\nu Y^{y,\nu}_t(\tilde y)\|^2\leq C_{0}\,\e^{-(c_2-c_1-\gamma)t},
\end{align*}
which in turn imply that
\begin{align*}
\mE\|Z^{\eta}_t(\tilde y)\|^2&=\mE\|\p_\nu Y^{y,\eta}_t(\tilde y)|_{y=\eta}\|^2=\mE\mE\big[\|\p_\nu Y^{y,\eta}_t(\tilde y)|_{y=\eta}\|^2|\cF_0\big]\\
&\leq\sup_{y\in\mR^{d_{2}}}\mE\|\p_\nu Y^{y,\nu}_t(\tilde y)\|^2\leq C_{0}\,\e^{-(c_2-c_1-\gamma)t},
\end{align*}
and
\begin{align*}
\mE\|\p_{\tilde y}Z^{\eta}_t(\tilde y)\|^2\leq C_{0}e^{-(c_2-c_1-\gamma)t}.
\end{align*}
Thus we arrive at
\begin{align*}
&\sup_{y\in\mR^{d_{2}},\nu\in\sP_{2}(\mR^{d_{2}}),\tilde y\in\mR^{d_{2}}}\Big(\|\p_\nu U(y,\nu)(\tilde y)\|+\|\p_{\tilde y}\p_\nu U(y,\nu)(\tilde y)\|\Big)\leq C_{0}<\infty.
\end{align*}
As a result, for every fixed $y\in\mR^{d_2}$ we obtain $U(y,\cdot)\in C_b^{(1,1)}(\sP_{2}(\mR^{d_{2}}))$.

 Similarly, we have
\begin{align*}
&\p_y\p_\nu U(y,\nu)(\tilde y)=\p_\nu\p_y U(y,\nu)(\tilde y)\\
&=\int_0^\infty\mE\Big[\p^2_yf(Y^{y,\nu}_t,\cL_{Y^\eta_t})\cdot\p_y Y^{y,\nu}_t\cdot\p_\nu Y^{y,\nu}_t(\tilde y)+\p_yf(Y^{y,\nu}_t,\cL_{Y^\eta_t})\cdot\p_\nu\p_y Y^{y,\nu}_t(\tilde y)\\
&\qquad\quad+\tilde\mE\big[\p_\nu\p_yf(Y^{y,\nu}_t,\cL_{Y^\eta_t})(\tilde Y^{\tilde y,\nu}_t)\cdot\p_y Y^{y,\nu}_t\cdot\p_y \tilde Y^{\tilde y,\nu}_t\big]\\
&\qquad\quad+\tilde\mE\big[\p_\nu\p_yf(Y^{y,\nu}_t,\cL_{Y^\eta_t})(\tilde Y^{\tilde \eta}_t)\cdot\p_y Y^{y,\nu}_t \cdot \tilde Z^{\tilde \eta}_t(\tilde y)\big]\Big]\dif t,
\end{align*}
where
\begin{align*}
\dif \p_\nu\p_y Y^{y,\nu}_t(\tilde y)&=\p^2_y b(Y^{y,\nu}_t,\cL_{Y^\eta_t})\cdot\p_y Y^{y,\nu}_t\cdot\p_\nu Y^{y,\nu}_t(\tilde y)\dif t+\p_y b(Y^{y,\nu}_t,\cL_{Y^\eta_t})\cdot\p_\nu\p_y Y^{y,\nu}_t(\tilde y)\dif t\\
&\quad+\tilde\mE\big[\p_\nu\p_y b(Y^{y,\nu}_t,\cL_{Y^\eta_t})(\tilde Y^{\tilde y,\nu}_t)\cdot\p_y Y_t^{y,\nu}\cdot\p_y \tilde Y^{\tilde y,\nu}_t\big]\dif t \\
&\quad+\tilde\mE\big[\p_\nu\p_y b(Y^{y,\nu}_t,\cL_{Y^\eta_t})(\tilde Y^{\tilde\eta}_t)\cdot\p_y Y_t^{y,\nu}\cdot \tilde Z^{\tilde \eta}_t(\tilde y)\big]\dif t\\
&\quad+\p^2_y \sigma(Y^{y,\nu}_t,\cL_{Y^\eta_t})\cdot\p_y Y^{y,\nu}_t\cdot\p_\nu Y^{y,\nu}_t(\tilde y) \dif W_t\\
&\quad+\p_y \sigma(Y^{y,\nu}_t,\cL_{Y^\eta_t})\cdot\p_\nu\p_y Y^{y,\nu}_t(\tilde y) \dif W_t\\
&\quad+\tilde\mE\big[\p_\nu\p_y \sigma(Y^{y,\nu}_t,\cL_{Y^\eta_t})(\tilde Y^{\tilde y,\nu}_t)\cdot\p_y Y_t^{y,\nu}\cdot\p_y \tilde Y^{\tilde y,\nu}_t\big]\dif W_t \\
&\quad+\tilde\mE\big[\p_\nu\p_y \sigma(Y^{y,\nu}_t,\cL_{Y^\eta_t})(\tilde Y^{\tilde\eta}_t)\cdot\p_y Y_t^{y,\nu}\cdot \tilde Z^{\tilde \eta}_t(\tilde y)\big]\dif W_t.
\end{align*}
By Lemma  \ref{Y2}  in the Appendix, we have
\begin{align*}
\sup_{y\in\mR^{d_{2}},\nu\in\sP_{2}(\mR^{d_{2}}),\tilde y\in\mR^{d_{2}}}\mE\|\p_\nu\p_y Y^{y,\nu}_t(\tilde y)\|^2\leq C_{0}\,\e^{-(c_2-c_1-\gamma)t},
\end{align*}
which in turn  implies
\begin{align*}
\sup_{y\in\mR^{d_{2}},\nu\in\sP_{2}(\mR^{d_{2}}),\tilde y\in\mR^{d_{2}}}\|\p_y\p_\nu U(y,\nu)(\tilde y)\|\leq C_{0}<\infty.
\end{align*}
Combing the above results, we obtain $U\in C_b^{2,(1,1)}(\mR^{d_{2}}\times\sP_{2}(\mR^{d_{2}}))$. The derivation that $U$ is the unique solution of equation (\ref{ellip2}) follows by It\^o's formula.

\vspace{1mm}
\noindent (iv) We proceed to prove the conclusions in (iii). The proof of $U(\cdot,\nu)\in C_b^{2k}(\mR^{d_2})$ when $\sigma, b, f\in C_b^{2k,(k,k)}(\mR^{d_2}\times\sP_2(\mR^{d_2}))$ is entirely similar  as in step (ii). Let us focus on the higher derivatives with respect to the $\nu$-variable. In view of (\ref{nuU}) and by the definition of the $L$-derivative, we have
\begin{align*}
&\p^2_\nu U(y,\nu)(\tilde y,\breve y)\\
&=\int_0^\infty\mE\Big[\p^2_yf(Y^{y,\nu}_t,\cL_{Y^\eta_t})\cdot\p_\nu Y^{y,\nu}_t(\tilde y)\cdot\p_\nu Y^{y,\nu}_t(\breve y)+\p_yf(Y^{y,\nu}_t,\cL_{Y^\eta_t})\cdot\p^2_\nu Y^{y,\nu}_t(\tilde y,\breve y)\\
&\qquad\quad+\breve\mE\big[\p_\nu\p_y f(Y^{y,\nu}_t,\cL_{Y^\eta_t})(\breve Y^{\breve y,\nu}_t)\cdot\p_\nu Y^{y,\nu}_t(\tilde y)\cdot\p_y \breve Y^{\breve y,\nu}_t\big] \\
&\qquad\quad+\breve\mE\big[\p_\nu\p_y f(Y^{y,\nu}_t,\cL_{Y^\eta_t})(\breve Y^{\breve\eta}_t)\cdot\p_\nu Y^{y,\nu}_t(\tilde y)\cdot\breve Z^{\breve\eta}_t(\breve y)\big] \\
&\qquad\quad+\tilde\mE\big[\p_\nu f(Y^{y,\nu}_t,\cL_{Y^\eta_t})(\tilde Y^{\tilde y,\nu}_t)\cdot\p_\nu\p_y \tilde Y^{\tilde y,\nu}_t(\breve y)\big] \\
&\qquad\quad+\tilde\mE\big[\p_{y}\p_\nu f(Y^{y,\nu}_t,\cL_{Y^\eta_t})(\tilde Y^{\tilde y,\nu}_t)\cdot\p_y \tilde Y^{\tilde y,\nu}_t\cdot\p_\nu Y^{y,\nu}_t(\breve y)\big] \\
&\qquad\quad+\tilde\mE\big[\p_{\tilde y}\p_\nu f(Y^{y,\nu}_t,\cL_{Y^\eta_t})(\tilde Y^{\tilde y,\nu}_t)\cdot\p_y \tilde Y^{\tilde y,\nu}_t\cdot\p_\nu \tilde Y^{y,\nu}_t(\breve y)\big] \\
&\qquad\quad+\tilde\mE\breve\mE\big[\p^2_\nu f(Y^{y,\nu}_t,\cL_{Y^\eta_t})(\tilde Y^{\tilde y,\nu}_t,\breve Y^{\breve y,\nu}_t)\cdot\p_y \tilde Y^{\tilde y,\nu}_t\cdot\p_y \breve Y^{\breve y,\nu}_t\big] \\
&\qquad\quad+\tilde\mE\breve\mE\big[\p^2_\nu f(Y^{y,\nu}_t,\cL_{Y^\eta_t})(\tilde Y^{\tilde y,\nu}_t,\breve Y^{\breve\eta}_t)\cdot\p_y \tilde Y^{\tilde y,\nu}_t\cdot\breve Z^{\breve\eta}_t(\breve y)\big]\\
&\qquad\quad+\tilde\mE\big[\p_\nu f(Y^{y,\nu}_t,\cL_{Y^\eta_t})(\tilde Y^{\breve y,\nu}_t)\cdot\p_y\p_\nu\tilde Y^{\breve y,\nu}_t(\tilde y)\big]+\tilde\mE\big[\p_\nu f(Y^{y,\nu}_t,\cL_{Y^\eta_t})(\tilde Y^{\tilde\eta}_t)\cdot\tilde Z^{\tilde\eta}_t(\tilde y,\breve y)\big]\\
&\qquad\quad+\tilde\mE\big[\p_{y}\p_\nu f(Y^{y,\nu}_t,\cL_{Y^\eta_t})(\tilde Y^{\tilde\eta}_t)\cdot\tilde Z^{\tilde\eta}_t(\tilde y)\cdot\p_\nu Y^{y,\nu}_t(\breve y)\big] \\
&\qquad\quad+\tilde\mE\big[\p_{\tilde y}\p_\nu f(Y^{y,\nu}_t,\cL_{Y^\eta_t})(\tilde Y^{\breve y,\nu}_t)\cdot\tilde Z^{\tilde\eta}_t(\tilde y)\cdot\p_y \tilde Y^{\breve y,\nu}_t\big] \\
&\qquad\quad+\tilde\mE\big[\p_{\tilde y}\p_\nu f(Y^{y,\nu}_t,\cL_{Y^\eta_t})(\tilde Y^{\tilde\eta}_t)\cdot\tilde Z^{\tilde\eta}_t(\tilde y)\cdot\tilde Z^{\tilde\eta}_t(\breve y)\big]\\
&\qquad\quad+\tilde\mE\breve\mE\big[\p^2_\nu f(Y^{y,\nu}_t,\cL_{Y^\eta_t})(\tilde Y^{\tilde\eta}_t,\breve Y^{\breve y,\nu}_t)\cdot \tilde Z^{ \tilde\eta}_t(\tilde y)\cdot\p_y \breve Y^{\breve y,\nu}_t\big] \\
&\qquad\quad+\tilde\mE\breve\mE\big[\p^2_\nu f(Y^{y,\nu}_t,\cL_{Y^\eta_t})(\tilde Y^{\tilde\eta}_t,\breve Y^{\breve\eta}_t)\cdot\tilde Z^{ \tilde\eta}_t(\tilde y)\cdot\breve Z^{\breve\eta}_t(\breve y)\big]\Big]dt,
\end{align*}
where $Z^\eta_t(\tilde y,\breve y ):=\p^2_\nu Y^{y,\nu}_t(\tilde y,\breve y)|_{y=\eta}$ with $\p^2_\nu Y^{y,\nu}_t(\tilde y,\breve y)$ satisfying
\begin{align*}
\dif\p^2_\nu Y^{y,\nu}_t(\tilde y,\breve y)&=\p^2_yb(Y^{y,\nu}_t,\cL_{Y^\eta_t})\cdot\p_\nu Y^{y,\nu}_t(\tilde y)\cdot\p_\nu Y^{y,\nu}_t(\breve y)\dif t\\
&\quad+\p_yb(Y^{y,\nu}_t,\cL_{Y^\eta_t})\cdot\p^2_\nu Y^{y,\nu}_t(\tilde y,\breve y)\dif t\\
&\quad+\p^2_y\sigma(Y^{y,\nu}_t,\cL_{Y^\eta_t})\cdot\p_\nu Y^{y,\nu}_t(\tilde y)\cdot\p_\nu Y^{y,\nu}_t(\breve y)\dif W_t\\
&\quad+\p_y\sigma(Y^{y,\nu}_t,\cL_{Y^\eta_t})\cdot\p^2_\nu Y^{y,\nu}_t(\tilde y,\breve y)\dif W_t+\sE_1(\tilde y, \breve y)+\sE_2(\tilde y, \breve y),
\end{align*}
and $\sE_1(\tilde y, \breve y)$ contains the drift part given by
\begin{align*}
\sE_1(\tilde y, \breve y):\!&=\breve\mE\big[\p_\nu\p_y b(Y^{y,\nu}_t,\cL_{Y^\eta_t})(\breve Y^{\breve y,\nu}_t)\cdot\p_\nu Y^{y,\nu}_t(\tilde y)\cdot\p_y \breve Y^{\breve y,\nu}_t\big]\dif t \\
&\quad+\breve\mE\big[\p_\nu\p_y b(Y^{y,\nu}_t,\cL_{Y^\eta_t})(\breve Y^{\breve\eta}_t)\cdot\p_\nu Y^{y,\nu}_t(\tilde y)\cdot\breve Z^{\breve\eta}_t(\breve y)\big]\dif t \\
&\quad+\tilde\mE\big[\p_\nu b(Y^{y,\nu}_t,\cL_{Y^\eta_t})(\tilde Y^{\tilde y,\nu}_t)\cdot\p_\nu\p_y \tilde Y^{\tilde y,\nu}_t(\breve y)\big]\dif t \\
&\quad+\tilde\mE\big[\p_{y}\p_\nu b(Y^{y,\nu}_t,\cL_{Y^\eta_t})(\tilde Y^{\tilde y,\nu}_t)\cdot\p_y \tilde Y^{\tilde y,\nu}_t\cdot\p_\nu Y^{y,\nu}_t(\breve y)\big]\dif t \\
&\quad+\tilde\mE\big[\p_{\tilde y}\p_\nu b(Y^{y,\nu}_t,\cL_{Y^\eta_t})(\tilde Y^{\tilde y,\nu}_t)\cdot\p_y \tilde Y^{\tilde y,\nu}_t\cdot\p_\nu \tilde Y^{y,\nu}_t(\breve y)\big] \dif t\\
&\quad+\tilde\mE\breve\mE\big[\p^2_\nu b(Y^{y,\nu}_t,\cL_{Y^\eta_t})(\tilde Y^{\tilde y,\nu}_t,\breve Y^{\breve y,\nu}_t)\cdot\p_y \tilde Y^{\tilde y,\nu}_t\cdot\p_y \breve Y^{\breve y,\nu}_t\big]\dif t \\
&\quad+\tilde\mE\breve\mE\big[\p^2_\nu b(Y^{y,\nu}_t,\cL_{Y^\eta_t})(\tilde Y^{\tilde y,\nu}_t,\breve Y^{\breve\eta}_t)\cdot\p_y \tilde Y^{\tilde y,\nu}_t\cdot\breve Z^{\breve\eta}_t(\breve y)\big]\dif t\\
&\quad+\tilde\mE\big[\p_\nu b(Y^{y,\nu}_t,\cL_{Y^\eta_t})(\tilde Y^{\breve y,\nu}_t)\cdot\p_y\p_\nu\tilde Y^{\breve y,\nu}_t(\tilde y)\big]\dif t\\
&\quad+\tilde\mE\big[\p_\nu b(Y^{y,\nu}_t,\cL_{Y^\eta_t})(\tilde Y^{\tilde\eta}_t)\cdot\tilde Z^{\tilde\eta}_t(\tilde y,\breve y)\big]\dif t\\
&\quad+\tilde\mE\big[\p_{y}\p_\nu b(Y^{y,\nu}_t,\cL_{Y^\eta_t})(\tilde Y^{\tilde\eta}_t)\cdot\tilde Z^{\tilde\eta}_t(\tilde y)\cdot\p_\nu Y^{y,\nu}_t(\breve y)\big] \dif t\\
&\quad+\tilde\mE\big[\p_{\tilde y}\p_\nu b(Y^{y,\nu}_t,\cL_{Y^\eta_t})(\tilde Y^{\breve y,\nu}_t)\cdot\tilde Z^{\tilde\eta}_t(\tilde y)\cdot\p_y \tilde Y^{\breve y,\nu}_t\big] \dif t\\
&\quad+\tilde\mE\big[\p_{\tilde y}\p_\nu b(Y^{y,\nu}_t,\cL_{Y^\eta_t})(\tilde Y^{\tilde\eta}_t)\cdot\tilde Z^{\tilde\eta}_t(\tilde y)\cdot\tilde Z^{\tilde\eta}_t(\breve y)\big]\dif t \\
&\quad+\tilde\mE\breve\mE\big[\p^2_\nu b(Y^{y,\nu}_t,\cL_{Y^\eta_t})(\tilde Y^{\tilde\eta}_t,\breve Y^{\breve y,\nu}_t)\cdot \tilde Z^{ \tilde\eta}_t(\tilde y)\cdot\p_y \breve Y^{\breve y,\nu}_t\big] \dif t\\
&\quad+\tilde\mE\breve\mE\big[\p^2_\nu b(Y^{y,\nu}_t,\cL_{Y^\eta_t})(\tilde Y^{\tilde\eta}_t,\breve Y^{\breve\eta}_t)\cdot\tilde Z^{ \tilde\eta}_t(\tilde y)\cdot\breve Z^{\breve\eta}_t(\breve y)\dif t,
\end{align*}
and $\sE_2(\tilde y, \breve y)$ involves the martingale terms given by
\begin{align*}
\sE_2(\tilde y, \breve y):\!&=\breve\mE\big[\p_\nu\p_y \sigma(Y^{y,\nu}_t,\cL_{Y^\eta_t})(\breve Y^{\breve y,\nu}_t)\cdot\p_\nu Y^{y,\nu}_t(\tilde y)\cdot\p_y \breve Y^{\breve y,\nu}_t\big]\dif W_t \\
&\quad+\breve\mE\big[\p_\nu\p_y \sigma(Y^{y,\nu}_t,\cL_{Y^\eta_t})(\breve Y^{\breve\eta}_t)\cdot\p_\nu Y^{y,\nu}_t(\tilde y)\cdot\breve Z^{\breve\eta}_t(\breve y)\big]\dif W_t \\
&\quad+\tilde\mE\big[\p_\nu \sigma(Y^{y,\nu}_t,\cL_{Y^\eta_t})(\tilde Y^{\tilde y,\nu}_t)\cdot\p_\nu\p_y \tilde Y^{\tilde y,\nu}_t(\breve y)\big]\dif W_t \\
&\quad+\tilde\mE\big[\p_{y}\p_\nu \sigma(Y^{y,\nu}_t,\cL_{Y^\eta_t})(\tilde Y^{\tilde y,\nu}_t)\cdot\p_y \tilde Y^{\tilde y,\nu}_t\cdot\p_\nu Y^{y,\nu}_t(\breve y)\big]\dif W_t \\
&\quad+\tilde\mE\big[\p_{\tilde y}\p_\nu \sigma(Y^{y,\nu}_t,\cL_{Y^\eta_t})(\tilde Y^{\tilde y,\nu}_t)\cdot\p_y \tilde Y^{\tilde y,\nu}_t\cdot\p_\nu \tilde Y^{y,\nu}_t(\breve y)\big] \dif W_t\\
&\quad+\tilde\mE\breve\mE\big[\p^2_\nu \sigma(Y^{y,\nu}_t,\cL_{Y^\eta_t})(\tilde Y^{\tilde y,\nu}_t,\breve Y^{\breve y,\nu}_t)\cdot\p_y \tilde Y^{\tilde y,\nu}_t\cdot\p_y \breve Y^{\breve y,\nu}_t\big]\dif W_t \\
&\quad+\tilde\mE\breve\mE\big[\p^2_\nu \sigma(Y^{y,\nu}_t,\cL_{Y^\eta_t})(\tilde Y^{\tilde y,\nu}_t,\breve Y^{\breve\eta}_t)\cdot\p_y \tilde Y^{\tilde y,\nu}_t\cdot\breve Z^{\breve\eta}_t(\breve y)\big]\dif W_t\\
&\quad+\tilde\mE\big[\p_\nu \sigma(Y^{y,\nu}_t,\cL_{Y^\eta_t})(\tilde Y^{\breve y,\nu}_t)\cdot\p_y\p_\nu\tilde Y^{\breve y,\nu}_t(\tilde y)\big]\dif W_t\\
&\quad+\tilde\mE\big[\p_\nu \sigma(Y^{y,\nu}_t,\cL_{Y^\eta_t})(\tilde Y^{\tilde\eta}_t)\cdot\tilde Z^{\tilde\eta}_t(\tilde y,\breve y)\big]\dif W_t\\
&\quad+\tilde\mE\big[\p_{y}\p_\nu \sigma(Y^{y,\nu}_t,\cL_{Y^\eta_t})(\tilde Y^{\tilde\eta}_t)\cdot\tilde Z^{\tilde\eta}_t(\tilde y)\cdot\p_\nu Y^{y,\nu}_t(\breve y)\big] \dif W_t\\
&\quad+\tilde\mE\big[\p_{\tilde y}\p_\nu \sigma(Y^{y,\nu}_t,\cL_{Y^\eta_t})(\tilde Y^{\breve y,\nu}_t)\cdot\tilde Z^{\tilde\eta}_t(\tilde y)\cdot\p_y \tilde Y^{\breve y,\nu}_t\big] \dif W_t\\
&\quad+\tilde\mE\big[\p_{\tilde y}\p_\nu \sigma(Y^{y,\nu}_t,\cL_{Y^\eta_t})(\tilde Y^{\tilde\eta}_t)\cdot\tilde Z^{\tilde\eta}_t(\tilde y)\cdot\tilde Z^{\tilde\eta}_t(\breve y)\big]\dif W_t \\
&\quad+\tilde\mE\breve\mE\big[\p^2_\nu \sigma(Y^{y,\nu}_t,\cL_{Y^\eta_t})(\tilde Y^{\tilde\eta}_t,\breve Y^{\breve y,\nu}_t)\cdot \tilde Z^{ \tilde\eta}_t(\tilde y)\cdot\p_y \breve Y^{\breve y,\nu}_t\big] \dif W_t\\
&\quad+\tilde\mE\breve\mE\big[\p^2_\nu \sigma(Y^{y,\nu}_t,\cL_{Y^\eta_t})(\tilde Y^{\tilde\eta}_t,\breve Y^{\breve\eta}_t)\cdot\tilde Z^{ \tilde\eta}_t(\tilde y)\cdot\breve Z^{\breve\eta}_t(\breve y)\dif W_t.
\end{align*}
Applying Lemma \ref{Y2}  in the Appendix, we have
\begin{align*}
&\sup_{y\in\mR^{d_{2}},\nu\in\sP_{2}(\mR^{d_{2}}),\tilde y\in\mR^{d_{2}},\breve y\in\mR^{d_{2}}}\mE\|\p^2_\nu Y^{y,\nu}_t(\tilde y,\breve y)\|^2\leq C_{0}\,\e^{-(c_2-c_1-\gamma)t},\\
&\sup_{\nu\in\sP_{2}(\mR^{d_{2}}),\tilde y\in\mR^{d_{2}},\breve y\in\mR^{d_{2}}}\mE\|Z^{\eta}_t(\tilde y,\breve y)\|^2\leq C_{0}\,\e^{-(c_2-c_1-\gamma)t},
\end{align*}
which in turn implies that
\begin{align*}
\sup_{y\in\mR^{d_{2}},\nu\in\sP_{2}(\mR^{d_{2}}),\tilde y\in\mR^{d_{2}},\breve y\in\mR^{d_{2}}}\|\p^2_\nu U(y,\nu)(\tilde y,\breve y)\|\leq C_{0}<\infty.
\end{align*}
The derivation that $(y,\tilde y,\breve y)\mapsto\p^2_\nu U(y,\nu)(\tilde y,\breve y)$ is in $C_b^2(\mR^{d_2}\times\mR^{d_2}\times\mR^{d_2})$ and $(y,\tilde y)\mapsto\p_\nu U(y,\nu)(\tilde y)$ is in $C_b^3(\mR^{d_2}\times\mR^{d_2})$ follows by the same argument as in steps (ii) and (iii). Thus we obtain $U\in C_b^{4,(2,2)}(\mR^{d_2}\times\mR^{d_2})$ when $\sigma, b, f\in C_b^{4,(2,2)}(\mR^{d_2}\times\sP_2(\mR^{d_2}))$. For general $k>2$, the proof follows by the same arguments, we omit the details here.
\end{proof}

\subsection{Poisson equations with parameters}

Now, we consider the following parameterized Poisson equation on the whole space $\mR^{d_2}\times\sP_2(\mR^{d_2})$:
\begin{align}\label{pof}
\sL_0U(x,\mu,y,\nu)=-f(x,\mu,y,\nu),
\end{align}
where $(x,\mu)\in\mR^{d_1}\times\sP_2(\mR^{d_1})$  are  parameters, and the operator $\sL_0$ is defined by (\ref{l0}). Recall that $\zeta^\mu$ is the unique invariant measure for the frozen McKean-Vlasov equation (\ref{sde1}).
As in the previous subsection, we assume that $f$ satisfies the following centering condition:
\begin{align}\label{cenf}
\int_{\mR^{d_2}}f(x,\mu,y,\zeta^{\mu})\zeta^{\mu}(\dif y)=0,\quad\forall (x,\mu)\in\mR^{d_1}\times\sP_2(\mR^{d_1}).
\end{align}
By regarding $(x,\mu)$ as parameters and according to Theorem \ref{pot}, the unique  solution $U$ to equation (\ref{pof})  should satisfy the centering condition (\ref{cenf}) and admit the probabilistic representation
\begin{align}\label{Pequ1}
U(x,\mu,y,\nu)&=\mE\left(\int_0^\infty f\big(x,\mu,Y_t^{\mu,y,\nu},\cL_{Y_t^{\mu,\eta}}\big)\dif t\right),
\end{align}
where $Y_t^{\mu,y,\nu}$ and $Y_t^{\mu,\eta}$ satisfy equation (\ref{sde1}) and (\ref{sde2}) with $\cL_\eta=\nu$, respectively.
The main problem addressed here is the smoothness of the solution $U$  with respect to the parameters $(x,\mu)$.
We have the following result, which will play a central role in our analysis below.

\bt\label{po}
Let  {\bf ($ H^{\sigma,b}$)} hold,  $\ell, j, k, m, n\in\mN$, the function $f$ satisfies the centering condition  (\ref{cenf}), and $U$ is defined by (\ref{Pequ1}).

\vspace{1mm}
(i) Assume that $a, b\in \mC_b^{(n,k),2(m-n),(m-n,m-n)}(\sP_2(\mR^{d_1})\times\mR^{d_2}\times\sP_2(\mR^{d_2}))$ and $f\in \mC_b^{j,(n,k),2(m-n),(m-n,m-n)}(\mR^{d_1}\times\sP_2(\mR^{d_1})\times\mR^{d_2}\times\sP_2(\mR^{d_2}))$ with $0\leq n\leq \ell<m$. Then we have
\begin{align}\label{U1}
U\in \mC_b^{j,(\ell,k),2(m-\ell),(m-\ell,m-\ell)}(\mR^{d_1}\times\sP_2(\mR^{d_1})\times\mR^{d_2}\times\sP_2(\mR^{d_2})).
\end{align}

\vspace{1mm}
(ii) Assume further that $a, b\in C_b^{(m,k),2m,(m,m)}(\sP_2(\mR^{d_1})\times\mR^{d_2}\times\sP_2(\mR^{d_2}))$ and  $f\in C_b^{j,(m,k),2m,(m,m)}(\mR^{d_1}\times\sP_2(\mR^{d_1})\times\mR^{d_2}\times\sP_2(\mR^{d_2}))$, then we    have
\begin{align}\label{U2}
U\in C_b^{j,(m,k),2m,(m,m)}(\mR^{d_1}\times\sP_2(\mR^{d_1})\times\mR^{d_2}\times\sP_2(\mR^{d_2})).
\end{align}

\et
\begin{proof}
	We devide the proof into three steps.
	
\vspace{1mm}	
\noindent (i) The regularities of $U$ with respect to the $(y,\nu)$-variables  follow  by Theorem \ref{pot}. Our task is to prove the regularities of $U$ with respect to $x$ and $\mu$.

The derivatives of $U$ with respect to the $x$-variable are easy. Since $f$ satisfies the centering condition (\ref{cenf}), it is obvious that $\p^j_xf$  satisfies (\ref{cenf}), too. Thus,  we can take derivative directly from both sides of (\ref{Pequ1}) to get that for any $j\in\mN$,
\begin{align}\label{Ux}
\p^j_xU(x,\mu,y,\nu)=\mE\left(\int_0^\infty \p^j_xf\big(x,\mu,Y_t^{\mu,y,\nu},\cL_{Y_t^{\mu,\eta}}\big)\dif t\right),
\end{align}
which in turn implies the desired conclusions for $U$ with respect to the $x$-variable.

\vspace{1mm}
(ii) The regularities of $U$ with respect to the $\mu$-variable is more delicate, since $\mu$ also appears in the process $Y_t^{\mu,y,\nu}$ as well as the distribution $\cL_{Y_t^{\mu,\eta}}$.
Taking derivative directly will involve complicated computations. We shall use the equation itself. Note that when $m=1$, the conclusion in (\ref{U1}) is obvious, we proceed to prove (\ref{U2}) with $m=1$.
Since $U$ is a classical solution to
$$
\sL_0(\mu,y,\nu)U(x,\mu,y,\nu)=-f(x,\mu,y,\nu),
$$
we have for every $\mu\in \sP_2(\mR^{d_1})$, $\phi\in L^2(\mR^{d_1})$ and $\rho>0$ that
\begin{align*}
&\sL_0(\mu,y,\nu)\Big(\frac{U(x,\mu,y,\nu)-U(x,\mu\circ(I+\rho\phi)^{-1},y,\nu)}{\rho}\Big)\\
&=\frac{f(x,\mu\circ(I+\rho\phi)^{-1},y,\nu)-f(x,\mu,y,\nu)}{\rho}\\
&\quad+\frac{1}{\rho}\Big(\sL_0(\mu\circ(I+\rho\phi)^{-1},y,\nu)-\sL_0(\mu,y,\nu)\Big)U(x,\mu\circ(I+\rho\phi)^{-1},y,\nu)\\
&=:h_1^{\phi}(x,\mu,y,\nu)(\rho),
\end{align*}
where
\begin{align*}
&\Big(\sL_0(\mu\circ(I+\rho\phi)^{-1},y,\nu)-\sL_0(\mu,y,\nu)\Big)U(x,\mu\circ(I+\rho\phi)^{-1},y,\nu)\\
&=\Big(b(\mu\circ(I+\rho\phi)^{-1},y,\nu)-b(\mu,y,\nu)\Big) \cdot\p_yU(x,\mu\circ(I+\rho\phi)^{-1},y,\nu)\\
&\quad+\frac{1}{2}\Big(a(\mu\circ(I+\rho\phi)^{-1},y,\nu)-a(\mu,y,\nu)\Big) \cdot\p^2_yU(x,\mu\circ(I+\rho\phi)^{-1},y,\nu)\\
&\quad+\int_{\mR^{d_2}}\Big[\big(b(\mu\circ(I+\rho\phi)^{-1},\tilde y,\nu)-b\big(\mu,\tilde y,\nu)\big)
\cdot\p_{\nu}U(x,\mu\circ(I+\rho\phi)^{-1},y,\nu)(\tilde y)\\
&\quad+\frac{1}{2}\big(a(\mu\circ(I+\rho\phi)^{-1},\tilde y,\nu)-a(\mu,\tilde y,\nu)\big)
\cdot\p_{\tilde y}\big[\p_{\nu}U(x,\mu\circ(I+\rho\phi)^{-1},y,\nu)(\tilde y)\big]\Big]\nu(\dif \tilde y).
\end{align*}
This implies that for every $\rho>0$, $h_1^{\phi}(x,\mu,y,\nu)(\rho)$ satisfies
$$
\int_{\mR^{d_2}}h_1^{\phi}(x,\mu,y,\zeta^\mu)(\rho)\zeta^\mu(\dif y)=0,
$$
and by Theorem \ref{pot} (i) we have the representation that
\begin{align}\label{d1}
&\frac{U(x,\mu\circ(I+\rho\phi)^{-1},y,\nu)-U(x,\mu,y,\nu)}{\rho}\no\\
&=\mE\left(\int_0^\infty h_1^{\phi}(x,\mu,Y_t^{\mu,y,\nu},\cL_{Y_t^{\mu,\eta}})(\rho)\dif t\right).
\end{align}
Note that by the definition of $L$-derivative,
\begin{align*}
&\lim_{\rho\to0}\frac{f(x,\mu\circ(I+\rho\phi)^{-1},y,\nu)-f(x,\mu,y,\nu)}{\rho}=\int_{\mR^{d_1}}\p_\mu f(x,\mu, y,\nu)(\tilde x)\cdot \phi(\tilde x)\mu(\dif \tilde x),\\
&\lim_{\rho\to0}\frac{b(\mu\circ(I+\rho\phi)^{-1},y,\nu)-b(\mu,y,\nu)}{\rho}=\int_{\mR^{d_1}}\p_\mu b(\mu,y,\nu)(\tilde x)\cdot \phi(\tilde x)\mu(\dif \tilde x),
\end{align*}
and
\begin{align*}
\lim_{\rho\to0}\frac{a(\mu\circ(I+\rho\phi)^{-1},y,\nu)-a(\mu,y,\nu)}{\rho}=\int_{\mR^{d_1}}\p_\mu a(\mu, y,\nu)(\tilde x)\cdot \phi(\tilde x)\mu(\dif \tilde x).
\end{align*}
By the continuity of $\p_y^2 U$ and $\p_{\tilde y}\p_\nu U$ we have
\begin{align*}
&\lim_{\rho\to0}h_1^\phi(x,\mu,y,\nu)(\rho)\\
&=\int_{\mR^{d_1}}\Big[\p_\mu f(x,\mu,y,\nu)(\tilde x)+\p_\mu b(\mu,y,\nu)(\tilde x)\cdot\p_yU(x,\mu,y,\nu)\no\\
&\qquad\qquad+\frac{1}{2}\p_\mu a(\mu,y,\nu)(\tilde x)\cdot\p^2_yU(x,\mu,y,\nu)\no\\
&\qquad\qquad+\int_{\mR^{d_2}}\Big[\p_\mu b(\mu,\tilde y,\nu)(\tilde x)\cdot\p_{\nu}U(x,\mu,y,\nu)(\tilde y)\no\\
&\qquad\qquad\qquad\quad+\frac{1}{2}\p_\mu a(\mu,\tilde y,\nu)(\tilde x)\cdot\p_{\tilde y}\big[\p_{\nu}U(x,\mu,y,\nu)(\tilde y)\big]\Big]\nu(\dif \tilde y)\Big]\cdot\phi(\tilde x)\mu(\dif\tilde x)\no\\
&=:\int_{\mR^{d_1}}h_1(x,\mu,y,\nu)(\tilde x)\cdot\phi(\tilde x)\mu(\dif\tilde x).
\end{align*}
This further implies that $h_1$ satisfies the centering condition (\ref{cenf}), i.e.,
\begin{align}\label{h1c}
\int_{\mR^{d_2}}h_1(x,\mu,y,\zeta^\mu)(\tilde x)\zeta^\mu(\dif y)=0.
\end{align}
Thus,  taking the limit $\rho\to0$ in (\ref{d1}),  we obtain
\begin{align*}
&\lim_{\rho\to0}\frac{U(x,\mu\circ(I+\rho\phi)^{-1},y,\nu)-U(x,\mu,y,\nu)}{\rho}\no\\
&=\mE\left(\int_0^\infty\int_{\mR^{d_1}}h_1(x,\mu,Y_t^{\mu,y,\nu},\cL_{Y_t^{\mu,\eta}})(\tilde x)\cdot\phi(\tilde x)\mu(\dif\tilde x)\dif t\right).
\end{align*}
As a result, $U(x,\mu,y,\nu)$ is $L$-differentiable at $\mu$ and by regarding $(x,\mu)$ as parameters, we arrive at
$$
\p_\mu U(x,\mu,y,\nu)(\tilde x)=\mE\left(\int_0^\infty h_1(x,\mu,Y_t^{\mu,y,\nu},\cL_{Y_t^{\mu,\eta}})(\tilde x)\dif t\right).
$$
Meanwhile, note that for every $k\in\mN$, $\p^k_{\tilde x} h_1$  also satisfies the centering condition, thus we have
\begin{align*}
\p^k_{\tilde x}\p_\mu U(x,\mu,y,\nu)(\tilde x)=\mE\left(\int_0^\infty \p^k_{\tilde x}h_1(x,\mu,Y_t^{\mu,y,\nu},\cL_{Y_t^{\mu,\eta}})(\tilde x)\dif t\right),
\end{align*}
where
\begin{align*}
\p^k_{\tilde x}h_1(x,\mu,y,\nu)(\tilde x)&=\p^k_{\tilde x}\p_\mu f(x,\mu,y,\nu)(\tilde x)+\p^k_{\tilde x}\p_\mu b(\mu,y,\nu)(\tilde x)\cdot\p_yU(x,\mu,y,\nu)\no\\
&\quad+\frac{1}{2}\p^k_{\tilde x}\p_\mu a(\mu,y,\nu)(\tilde x)\cdot\p^2_yU(x,\mu,y,\nu)\no\\
&\quad+\int_{\mR^{d_2}}\Big[\p^k_{\tilde x}\p_\mu b(\mu,\tilde y,\nu)(\tilde x)\cdot\p_{\nu}U(x,\mu,y,\nu)(\tilde y)\no\\
&\qquad\quad+\frac{1}{2}\p^k_{\tilde x}\p_\mu a(\mu,\tilde y,\nu)(\tilde x)\cdot\p_{\tilde y}\big[\p_{\nu}U(x,\mu,y,\nu)(\tilde y)\big]\Big]\nu(\dif \tilde y).
\end{align*}
Consequently, we have $U(x,\cdot,y,\nu)\in C_b^{(1,k)}(\sP_2(\mR^{d_1}))$.

\vspace{1mm}
(iii) Now we prove (\ref{U1}) and (\ref{U2}) with $m=2$. By the assumptions that $a,b\in \mC_b^{(1,k),2,(1,1)}(\sP_2(\mR^{d_1})\times\mR^{d_2}\times\sP_2(\mR^{d_2}))$ and $f\in \mC_b^{j,(1,k),2,(1,1)}(\mR^{d_1}\times\sP_2(\mR^{d_1})\times\mR^{d_2}\times\sP_2(\mR^{d_2}))$, we have
$$
\p^k_{\tilde x}\p_\mu a(\mu,\cdot,\cdot)(\tilde x), \p^k_{\tilde x}\p_\mu b(\mu,\cdot,\cdot)(\tilde x),\p^k_{\tilde x}\p_\mu f(x,\mu,\cdot,\cdot)(\tilde x)\in C_b^{2,(1,1)}(\mR^{d_2}\times\sP_2(\mR^{d_2})).
 $$
 Combing this with the fact that  $\p_y^2U(x,\mu,\cdot,\cdot), \p_{\tilde y}\p_\nu U(x,\mu,\cdot,\cdot)(\tilde y)\in C_b^{2,(1,1)}(\mR^{d_2}\times\sP_2(\mR^{d_2}))$, one can check that for every fixed $x,\tilde x\in\mR^{d_1}$ and $\mu\in\sP_2(\mR^{d_1})$, we have
$$
\p^k_{\tilde x}h_1(x,\mu,\cdot,\cdot)(\tilde x)\in C_b^{2,(1,1)}(\mR^{d_2}\times\sP_2(\mR^{d_2})).
$$
As a result of Theorem \ref{pot}, we conclude that
\begin{align*}
\sL_0(\mu,y,\nu)\p_\mu U(x,\mu,y,\nu)(\tilde x)=-h_1(x,\mu,y,\nu)(\tilde x)
\end{align*}
and
\begin{align*}
\sL_0(\mu,y,\nu)\p^k_{\tilde x}\p_\mu U(x,\mu,y,\nu)(\tilde x)=-\p_{\tilde x}h_1(x,\mu,y,\nu)(\tilde x),
\end{align*}
which in turn yields that $\p^k_{\tilde x}\p_\mu U(x,\mu,\cdot,\cdot)(\tilde x)\in C_b^{2,(1,1)}(\mR^{d_2}\times\sP_2(\mR^{d_2}))$. Consequently, we get $U(x,\cdot,\cdot,\cdot)\in\mC_b^{(1,k),2,(1,1)}(\sP_2(\mR^{d_1})\times\mR^{d_2}\times\sP_2(\mR^{d_2}))$. In view of (\ref{Ux}), by $\p_x^j f(x,\cdot,\cdot,\cdot)\in \mC_b^{(1,k),2,(1,1)}(\sP_2(\mR^{d_1})\times\mR^{d_2}\times\sP_2(\mR^{d_2}))$ and Theorem \ref{pot}, we obtain
\begin{align*}
\sL_0(\mu,y,\nu)\p_x^j U(x,\mu,y,\nu)=-\p_x^jf(x,\mu,y,\nu),
\end{align*}
which in turns implies that $\p_x^j U(x,\cdot,\cdot,\cdot)\in \mC_b^{(1,k),2,(1,1)}(\sP_2(\mR^{d_1})\times\mR^{d_2}\times\sP_2(\mR^{d_2}))$. Thus we have
 $$
 U\in \mC_b^{j,(1,k),2,(1,1)}(\mR^{d_1}\times\sP_2(\mR^{d_1})\times\mR^{d_2}\times\sP_2(\mR^{d_2})).
 $$
It remains to show the second order derivatives of  $U(x,\mu,y,\nu)$ with respect to the $\mu$-variable.
The same argument as above we have
$$
\p^2_\mu U(x,\mu,y,\nu)(\tilde x,\tilde{\tilde x})=\mE\left(\int_0^\infty h_2(x,\mu,Y_t^{\mu,y,\nu},\cL_{Y_t^{\mu,\eta}})(\tilde x,\tilde{\tilde x})\dif t\right),
$$
where $h_2$ satisfies the centering condition (\ref{cenf}) and is given by
\begin{align*}
&h_2(x,\mu,y,\nu)(\tilde x,\tilde{\tilde x}):=\p^2_\mu f(x,\mu,y,\nu)(\tilde x,\tilde{\tilde x})+\p^2_\mu b(\mu,y,\nu)(\tilde x,\tilde{\tilde x})\cdot\p_y U(x,\mu,y,\nu)\\
&\quad+\p_\mu b(\mu,y,\nu)(\tilde x)\cdot\p_\mu\p_y U(x,\mu,y,\nu)(\tilde{\tilde x})+\p_\mu b(\mu,y,\nu)(\tilde{\tilde x})\cdot\p_y [\p_\mu U(x,\mu,y,\nu)(\tilde x)]\\
&\quad+\frac{1}{2}\p_\mu a(\mu,y,\nu)(\tilde{\tilde x})\cdot\p^2_y[\p_\mu U(x,\mu,y,\nu)(\tilde x)]\\
&\quad+\frac{1}{2}\p^2_\mu a(\mu,y,\nu)(\tilde x,\tilde{\tilde x})\cdot\p^2_y U(x,\mu,y,\nu)+\frac{1}{2}\p_\mu a(\mu,y,\nu)(\tilde x)\cdot\p_\mu\p^2_y U(x,\mu,y,\nu)(\tilde{\tilde x})\\
&\quad+\int_{\mR^{d_2}}\Big[\p_\mu b(\mu,\tilde y,\nu)(\tilde{\tilde x})\cdot\p_{\nu}[\p_\mu U(x,\mu,y,\nu)(\tilde x)](\tilde y)\\
&\qquad\qquad+\p^2_\mu b(\mu,\tilde y,\nu)(\tilde x,\tilde{\tilde x})\cdot\p_{\nu}U(x,\mu,y,\nu)(\tilde y)\\
&\qquad\qquad+\p_\mu b(\mu,\tilde y,\nu)(\tilde x)\cdot\p_{\mu}[\p_\nu U(x,\mu,y,\nu)(\tilde y)](\tilde{\tilde x})\Big]\nu(\dif \tilde y)\\
&\quad+\frac{1}{2}\int_{\mR^{d_2}}\!\!\Big[\p_\mu a(\mu,\tilde y,\nu)(\tilde{\tilde x})\cdot\p_{\tilde y}\big[\p_{\nu}[\p_\mu U(x,\mu,y,\nu)(\tilde x)](\tilde y)\big]\\
&\qquad\qquad\quad+\p^2_\mu a(\mu,\tilde y,\nu)(\tilde x,\tilde{\tilde x})\cdot\p_{\tilde y}[\p_{\nu}U(x,\mu,y,\nu)(\tilde y)]\\
&\qquad\qquad\quad+\p_\mu a(\mu,\tilde y,\nu)(\tilde x)\cdot\p_{\mu}\big[\p_{\tilde y}[\p_\nu U(x,\mu,y,\nu)(\tilde y)](\tilde{\tilde x})\big]\Big]\nu(\dif \tilde y).
\end{align*}
Consequently, we deduce that
$$
\p_{\tilde{\tilde x}}[\p^2_\mu U(x,\mu,y,\nu)(\tilde x,\tilde{\tilde x})]=\mE\left(\int_0^\infty \p_{\tilde{\tilde x}}h_2(x,\mu,Y_t^{\mu,y,\nu},\cL_{Y_t^{\mu,\eta}})(\tilde x,\tilde{\tilde x})\dif t\right),
$$
where
\begin{align*}
\p_{\tilde{\tilde x}}h_2(x,\mu,y,\nu)(\tilde x,\tilde{\tilde x})&=\p_{\tilde{\tilde x}}[\p^2_\mu f(x,\mu,y,\nu)(\tilde x,\tilde{\tilde x})]+\p_{\tilde{\tilde x}}[\p^2_\mu b(\mu,y,\nu)(\tilde x,\tilde{\tilde x})]\cdot\p_y U(x,\mu,y,\nu)\\
&\quad+\p_\mu b(\mu,y,\nu)(\tilde x)\cdot\p_{\tilde{\tilde x}}[\p_\mu\p_y U(x,\mu,y,\nu)(\tilde{\tilde x})]\\
&\quad+\p_{\tilde{\tilde x}}[\p_\mu b(\mu,y,\nu)(\tilde{\tilde x})]\cdot\p_y [\p_\mu U(x,\mu,y,\nu)(\tilde x)]\\
&\quad+\frac{1}{2}\p_{\tilde{\tilde x}}[\p_\mu a(\mu,y,\nu)(\tilde{\tilde x})]\cdot\p^2_y[\p_\mu U(x,\mu,y,\nu)(\tilde x)]\\
&\quad+\frac{1}{2}\p_{\tilde{\tilde x}}[\p^2_\mu a(\mu,y,\nu)(\tilde x,\tilde{\tilde x})]\cdot\p^2_y U(x,\mu,y,\nu)\\
&\quad+\frac{1}{2}\p_\mu a(\mu,y,\nu)(\tilde x)\cdot\p_{\tilde{\tilde x}}[\p_\mu\p^2_y U(x,\mu,y,\nu)(\tilde{\tilde x})]\\
&\quad+\int_{\mR^{d_2}}\Big[\p_{\tilde{\tilde x}}[\p_\mu b(\mu,\tilde y,\nu)(\tilde{\tilde x})]\cdot\p_{\nu}[\p_\mu U(x,\mu,y,\nu)(\tilde x)](\tilde y)\\
&\qquad\qquad+\p_{\tilde{\tilde x}}[\p^2_\mu b(\mu,\tilde y,\nu)(\tilde x,\tilde{\tilde x})]\cdot\p_{\nu}U(x,\mu,y,\nu)(\tilde y)\\
&\qquad\qquad+\p_\mu b(\mu,\tilde y,\nu)(\tilde x)\cdot\p_{\tilde{\tilde x}}\big[\p_{\mu}[\p_\nu U(x,\mu,y,\nu)(\tilde y)](\tilde{\tilde x})\big]\Big]\nu(\dif \tilde y)\\
&\quad+\frac{1}{2}\int_{\mR^{d_2}}\Big[\p_{\tilde{\tilde x}}[\p_\mu a(\mu,\tilde y,\nu)(\tilde{\tilde x})]\cdot\p_{\tilde y}\big[\p_{\nu}[\p_\mu U(x,\mu,y,\nu)(\tilde x)](\tilde y)\big]\\
&\qquad\qquad\quad+\p_{\tilde{\tilde x}}[\p^2_\mu a(\mu,\tilde y,\nu)(\tilde x,\tilde{\tilde x})]\cdot\p_{\tilde y}[\p_{\nu}U(x,\mu,y,\nu)(\tilde y)]\\
&\qquad\qquad\quad+\p_\mu a(\mu,\tilde y,\nu)(\tilde x)\cdot\p_{\tilde{\tilde x}}\big[\p_{\mu}\p_{\tilde y}[\p_\nu U(x,\mu,y,\nu)(\tilde y)](\tilde{\tilde x})\big]\Big]\nu(\dif \tilde y).
\end{align*}
As a result, we obtain $U(x,\cdot,y,\nu)\in C_b^{(2,1)}(\sP_2(\mR^{d_1}))$. In the same way, we get $U(x,\cdot,y,\nu)\in C_b^{(2,k)}(\sP_2(\mR^{d_1}))$. For general $m\geq 3$, the proof follows by the similar argumants, we omit the details here.
\end{proof}

Given a function $f(x,\mu,y,\nu)$, we denote
\begin{align}\label{barf}
\bar f(x,\mu):=\int_{\mR^{d_2}}f(x,\mu,y,\zeta^\mu)\zeta^\mu(\dif y),
\end{align}
where $\zeta^\mu$ is the unique invariant measure for the frozen McKean-Vlasov equation (\ref{sde1}). The verification that the averaged coefficient $\bar f(x,\mu)$ are smooth usually constitutes a separate problem connected with the  smoothness of the invariant measure $\zeta^\mu$ with respect to the parameter $\mu$. Here,
as a direct application of Theorem \ref{po}, we have the following result.

\bc\label{avef}
Assume that  {\bf ($ H^{\sigma,b}$)} holds and  $\ell, m, k\in\mN$. If for every $1\leq n<m$, $a, b\in (C_b^{(m,k),2m,(m,m)}\cap \mC_b^{(n,k),2(m-n),(m-n,m-n)})(\sP_2(\mR^{d_1})\times\mR^{d_2}\times\sP_2(\mR^{d_2}))$ and  $f\in (C_b^{\ell,(m,k),2m,(m,m)}\cap \mC_b^{\ell,(n,k),2(m-n),(m-n,m-n)})(\mR^{d_1}\times\sP_2(\mR^{d_1})\times\mR^{d_2}\times\sP_2(\mR^{d_2}))$. Then we have
$\bar f\in C_b^{\ell,(m,k)}(\mR^{d_1}\times\sP_2(\mR^{d_1}))$.

In particular, we have
\begin{align}\label{barf1}
\p_\mu \bar f(x,\mu)(\tilde x)&=\int_{\mR^{d_2}}\bigg[\p_\mu f(x,\mu,y,\zeta^\mu)(\tilde x)+\p_\mu b(\mu,y,\zeta^\mu)(\tilde x)\cdot\p_yU(x,\mu,y,\zeta^\mu)\no\\
&\qquad+\frac{1}{2}\p_\mu a(\mu,y,\zeta^\mu)(\tilde x)\cdot\p^2_yU(x,\mu,y,\zeta^\mu)\no\\
&\qquad+\int_{\mR^{d_2}}\Big[\p_\mu b(\mu,\tilde y,\zeta^\mu)(\tilde x)\cdot\p_{\nu}U(x,\mu,y,\zeta^\mu)(\tilde y)\no\\
&\qquad +\frac{1}{2}\p_\mu a(\mu,\tilde y,\nu)(\tilde x)\cdot\p_{\tilde y}\big[\p_{\nu}U(x,\mu,y,\zeta^\mu)(\tilde y)\big]\Big]\zeta^\mu(\dif \tilde y)\bigg]\zeta^\mu(\dif y),
\end{align}
and
\begin{align}\label{barf2}
&\p^2_\mu\bar f(x,\mu)(\tilde x,\tilde{\tilde x})\no\\
&=\int_{\mR^{d_2}}\bigg[\p^2_\mu f(x,\mu,y,\zeta^\mu)(\tilde x,\tilde{\tilde x})+\p^2_\mu b(\mu,y,\zeta^\mu)(\tilde x,\tilde{\tilde x})\cdot\p_y U(x,\mu,y,\zeta^\mu)\no\\
&\qquad+\p_\mu b(\mu,y,\zeta^\mu)(\tilde x)\cdot\p_\mu\p_y U(x,\mu,y,\zeta^\mu)(\tilde{\tilde x})+\p_\mu b(\mu,y,\zeta^\mu)(\tilde{\tilde x})\cdot\p_y [\p_\mu U(x,\mu,y,\zeta^\mu)(\tilde x)]\no\\
&\qquad+\frac{1}{2}\p_\mu a(\mu,y,\zeta^\mu)(\tilde{\tilde x})\cdot\p^2_y[\p_\mu U(x,\mu,y,\zeta^\mu)(\tilde x)]+\frac{1}{2}\p^2_\mu a(\mu,y,\zeta^\mu)(\tilde x,\tilde{\tilde x})\cdot\p^2_y U(x,\mu,y,\zeta^\mu)\no\\
&\qquad+\frac{1}{2}\p_\mu a(\mu,y,\zeta^\mu)(\tilde x)\cdot\p_\mu\p^2_y U(x,\mu,y,\zeta^\mu)(\tilde{\tilde x})\no\\
&\qquad+\int_{\mR^{d_2}}\Big[\p_\mu b(\mu,\tilde y,\zeta^\mu)(\tilde{\tilde x})\cdot\p_{\nu}[\p_\mu U(x,\mu,y,\zeta^\mu)(\tilde x)](\tilde y)\no\\
&\qquad\qquad\quad+\p^2_\mu b(\mu,\tilde y,\zeta^\mu)(\tilde x,\tilde{\tilde x})\cdot\p_{\nu}U(x,\mu,y,\zeta^\mu)(\tilde y)\no\\
&\qquad\qquad\quad+\p_\mu b(\mu,\tilde y,\zeta^\mu)(\tilde x)\cdot\p_{\mu}[\p_\nu U(x,\mu,y,\zeta^\mu)(\tilde y)](\tilde{\tilde x})\Big]\zeta^\mu(\dif \tilde y)\no\\
&\qquad+\frac{1}{2}\int_{\mR^{d_2}}\Big[\p_\mu a(\mu,\tilde y,\zeta^\mu)(\tilde{\tilde x})\cdot\p_{\tilde y}\big[\p_{\nu}[\p_\mu U(x,\mu,y,\zeta^\mu)(\tilde x)](\tilde y)\big]\no\\
&\qquad\qquad\quad+\p^2_\mu a(\mu,\tilde y,\zeta^\mu)(\tilde x,\tilde{\tilde x})\cdot\p_{\tilde y}[\p_{\nu}U(x,\mu,y,\zeta^\mu)(\tilde y)]\no\\
&\qquad\qquad\quad+\p_\mu a(\mu,\tilde y,\zeta^\mu)(\tilde x)\cdot\p_{\mu}\p_{\tilde y}[\p_\nu U(x,\mu,y,\zeta^\mu)(\tilde y)](\tilde{\tilde x})\Big]\zeta^\mu(\dif \tilde y)\bigg]\zeta^\mu(\dif y),
\end{align}
where $U(x,\mu,y,\nu)$ is the unique solution to the Poisson equation
\begin{align}\label{Uf}
\sL_0U(x,\mu,y,\nu)=-\big[f(x,\mu,y,\nu)-\bar f(x,\mu)\big].
\end{align}
\ec
\begin{proof}
Since $\zeta^\mu$ does not depend on the $x$-variable,
we can take derivative directly with respect to  $x$ from both sides of (\ref{barf}) to get
\begin{align*}
&\p^\ell_x\bar f(x,\mu)=\int_{\mR^{d_2}}\p^\ell_xf(x,\mu,y,\zeta^\mu)\zeta^\mu(\dif y),
\end{align*}
which in turn implies that $\bar f(\cdot,\mu)\in C_b^\ell(\mR^{d_1})$. We proceed to show the regularities of $\bar f$ with respect to the $\mu$-variable.
Note that the function
$$
\delta f(x,\mu,y,\nu):=f(x,\mu,y,\nu)-\bar f(x,\mu)
$$
always satisfies the centering condition. Thus, under our assumptions there exists a unique solution $U$ to the Poisson equation (\ref{Uf}).
Following the same arguments as in (\ref{h1c}) we have
\begin{align}\label{fc}
&\int_{\mR^{d_2}}\bigg[\p_\mu \delta f(x,\mu,y,\zeta^\mu)(\tilde x)+\p_\mu b(\mu,y,\zeta^\mu)(\tilde x)\cdot\p_yU(x,\mu,y,\zeta^\mu)\no\\
&\qquad+\frac{1}{2}\p_\mu a(\mu,y,\zeta^\mu)(\tilde x)\cdot\p^2_yU(x,\mu,y,\zeta^\mu)\no\\
&\qquad+\int_{\mR^{d_2}}\Big[\p_\mu b(\mu,\tilde y,\zeta^\mu)(\tilde x)\cdot\p_{\nu}U(x,\mu,y,\zeta^\mu)(\tilde y)\no\\
&\qquad\qquad+\frac{1}{2}\p_\mu a(\mu,\tilde y,\nu)\cdot\p_{\tilde y}\big[\p_{\nu}U(x,\mu,y,\zeta^\mu)(\tilde y)\big]\Big]\zeta^\mu(\dif \tilde y)\bigg]\zeta^\mu(\dif y)=0.
\end{align}
Since
\begin{align*}
&\int_{\mR^{d_2}}\p_\mu \delta f(x,\mu,y,\zeta^\mu)(\tilde x)\zeta^\mu(\dif y)\\
&=-\p_\mu \bar f(x,\mu)(\tilde x)+\int_{\mR^{d_2}}\p_\mu f(x,\mu,y,\zeta^\mu)(\tilde x)\zeta^\mu(\dif y),
\end{align*}
this together with (\ref{fc}) yields (\ref{barf1}).
Taking derivatives with respect to the $\tilde x$ and $x$ from both sides of (\ref{barf1}), we obtain
\begin{align*}
\p_{\tilde x}[\p_\mu \bar f(x,\mu)(\tilde x)]&=\int_{\mR^{d_2}}\bigg[\p_{\tilde x}[\p_\mu f(x,\mu,y,\zeta^\mu)(\tilde x)]+\p_{\tilde x}[\p_\mu b(\mu,y,\zeta^\mu)(\tilde x)]\cdot\p_yU(x,\mu,y,\zeta^\mu)\\
&\qquad+\frac{1}{2}\p_{\tilde x}[\p_\mu a(\mu,y,\zeta^\mu)(\tilde x)]\cdot\p^2_yU(x,\mu,y,\zeta^\mu)\\
&\qquad+\int_{\mR^{d_2}}\Big[\p_{\tilde x}[\p_\mu b(\mu,\tilde y,\zeta^\mu)(\tilde x)]\cdot\p_{\nu}U(x,\mu,y,\zeta^\mu)(\tilde y)\\
&\qquad\qquad+\frac{1}{2}\p_{\tilde x}[\p_\mu a(\mu,\tilde y,\nu)(\tilde x)]\cdot\p_{\tilde y}\big[\p_{\nu}U(x,\mu,y,\zeta^\mu)(\tilde y)\big]\Big]\zeta^\mu(\dif \tilde y)\bigg]\zeta^\mu(\dif y),
\end{align*}
and
\begin{align*}
\p_x[\p_\mu \bar f(x,\mu)(\tilde x)]&=\int_{\mR^{d_2}}\bigg[\p_x[\p_\mu f(x,\mu,y,\zeta^\mu)(\tilde x)]+\p_\mu b(\mu,y,\zeta^\mu)(\tilde x)\cdot\p_x\p_yU(x,\mu,y,\zeta^\mu)\\
&\qquad+\frac{1}{2}\p_\mu a(\mu,y,\zeta^\mu)(\tilde x)\cdot\p_x\p^2_yU(x,\mu,y,\zeta^\mu)\\
&\qquad+\int_{\mR^{d_2}}\Big[\p_\mu b(\mu,\tilde y,\zeta^\mu)(\tilde x)\cdot\p_x[\p_{\nu}U(x,\mu,y,\zeta^\mu)(\tilde y)]\\
&\qquad\qquad+\frac{1}{2}\p_\mu a(\mu,\tilde y,\nu)(\tilde x)\cdot\p_{\tilde y}\big[\p_x[\p_{\nu}U(x,\mu,y,\zeta^\mu)(\tilde y)]\big]\Big]\zeta^\mu(\dif \tilde y)\bigg]\zeta^\mu(\dif y).
\end{align*}
Therefore, we deduce $\bar f\in C_b^{2,(1,1)}(\mR^{d_1}\times\sP_2(\mR^{d_1}))$.
Similarly, we have
\begin{align*}
&\int_{\mR^{d_2}}\bigg[\p^2_\mu \delta f(x,\mu,y,\zeta^\mu)(\tilde x,\tilde{\tilde x})+\p^2_\mu b(\mu,y,\zeta^\mu)(\tilde x,\tilde{\tilde x})\cdot\p_y U(x,\mu,y,\zeta^\mu)\\
&\qquad+\p_\mu b(\mu,y,\zeta^\mu)(\tilde x)\cdot\p_\mu\p_y U(x,\mu,y,\zeta^\mu)(\tilde{\tilde x})+\p_\mu b(\mu,y,\zeta^\mu)(\tilde{\tilde x})\cdot\p_y [\p_\mu U(x,\mu,y,\zeta^\mu)(\tilde x)]\\
&\qquad+\frac{1}{2}\p_\mu a(\mu,y,\zeta^\mu)(\tilde{\tilde x})\cdot\p^2_y[\p_\mu U(x,\mu,y,\zeta^\mu)(\tilde x)]+\frac{1}{2}\p^2_\mu a(\mu,y,\zeta^\mu)(\tilde x,\tilde{\tilde x})\cdot\p^2_y U(x,\mu,y,\zeta^\mu)\\
&\qquad+\frac{1}{2}\p_\mu a(\mu,y,\zeta^\mu)(\tilde x)\cdot\p_\mu\p^2_y U(x,\mu,y,\zeta^\mu)(\tilde{\tilde x})\\
&\qquad+\int_{\mR^{d_2}}\Big[\p_\mu b(\mu,\tilde y,\zeta^\mu)(\tilde{\tilde x})\cdot\p_{\nu}[\p_\mu U(x,\mu,y,\zeta^\mu)(\tilde x)](\tilde y)\\
&\qquad\qquad\quad+\p^2_\mu b(\mu,\tilde y,\zeta^\mu)(\tilde x,\tilde{\tilde x})\cdot\p_{\nu}U(x,\mu,y,\zeta^\mu)(\tilde y)\\
&\qquad\qquad\quad+\p_\mu b(\mu,\tilde y,\zeta^\mu)(\tilde x)\cdot\p_{\mu}[\p_\nu U(x,\mu,y,\zeta^\mu)(\tilde y)](\tilde{\tilde x})\Big]\zeta^\mu(\dif \tilde y)\\
&\qquad+\frac{1}{2}\int_{\mR^{d_2}}\Big[\p_\mu a(\mu,\tilde y,\zeta^\mu)(\tilde{\tilde x})\cdot\p_{\tilde y}\big[\p_{\nu}[\p_\mu U(x,\mu,y,\zeta^\mu)(\tilde x)](\tilde y)\big]\\
&\qquad\qquad\quad+\p^2_\mu a(\mu,\tilde y,\zeta^\mu)(\tilde x,\tilde{\tilde x})\cdot\p_{\tilde y}[\p_{\nu}U(x,\mu,y,\zeta^\mu)(\tilde y)]\\
&\qquad\qquad\quad+\p_\mu a(\mu,\tilde y,\zeta^\mu)(\tilde x)\cdot\p_{\mu}\p_{\tilde y}[\p_\nu U(x,\mu,y,\zeta^\mu)(\tilde y)](\tilde{\tilde x})\Big]\zeta^\mu(\dif \tilde y)\bigg]\zeta^\mu(\dif y)=0,
\end{align*}
which implies that (\ref{barf2}) holds. For general $m\in\mN$, the proof follows by induction, we omit the details here.
\end{proof}

\section{Fluctuations estimates}

Let $(X_t^\eps, Y_t^\eps)$ and $\bar X_t$ satisfy the McKean-Vlasov equations (\ref{sde}) and (\ref{ave}), respectively. This section is devoted to establish some integral estimates for the fluctuations between $X_t^\eps$ and its homogenized limit $\bar X_t$.
Given a function $f(t,x,\mu,y,\nu)$ on $[0,\infty)\times\mR^{d_1}\times\sP_2(\mR^{d_1})\times\mR^{d_2}\times\sP_2(\mR^{d_2})$, we consider the Poisson equation:
\begin{align}\label{pss}
\sL_0\psi(t,x,\mu,y,\nu)=-f(t,x,\mu,y,\nu),
\end{align}
where the operator $\sL_0$ is defined by (\ref{l0}), and $(t,x,\mu)\in \mR_+\times\mR^{d_1}\times\sP_2(\mR^{d_1})$ are regarded as parameters.
We first give the following  estimate of the functional law of large number type.

\bl\label{flu1}
Assume that {\bf (H$^{\sigma,b}$)} hold, $F,H,G,c\in C_b^{2,(1,1),2,(1,1)}(\mR^{d_1}\times\sP_2(\mR^{d_1})\times\mR^{d_2}\times\sP_2(\mR^{d_2}))$ and $b,\sigma\in C_b^{(1,1),2,(1,1)}(\sP_2(\mR^{d_1})\times\mR^{d_2}\times\sP_2(\mR^{d_2}))$. Then for every $f\in C_b^{1,2,(1,1),2,(1,1)}(\mR_+\times\mR^{d_1}\times\sP_2(\mR^{d_1})\times\mR^{d_2}\times\sP_2(\mR^{d_2}))$ satisfying (\ref{cen}), we have
\begin{align*}
\left|\mE\left(\int_0^tf(s,X_s^\eps,\cL_{X_s^\eps},Y_s^\eps,\cL_{Y_s^\eps})\dif s\right)\right|\leq C_0\,\eps,
\end{align*}
where $C_0>0$ is a constant independent of $\eps$.
\el
\begin{proof}
Let $\psi$ be the solution of the equation (\ref{pss}). According the assumptions on the coefficients  and Theorem \ref{po}, we have
 $\psi\in C_b^{1,2,(1,1),2,(1,1)}(\mR_+\times\mR^{d_1}\times\sP_2(\mR^{d_1})\times\mR^{d_2}\times\sP_2(\mR^{d_2}))$. Using It\^o's formula, we deduce that
\begin{align}
&\psi(t,X_t^\eps,\cL_{X_t^\eps},Y_t^\eps,\cL_{Y_t^\eps})\no\\
&=\psi(0,\xi,\cL_\xi,\eta,\cL_\eta) +\int_0^t\Big[\p_s\psi(s,X_s^\eps,\cL_{X_s^\eps},Y_s^\eps,\cL_{Y_s^\eps})\no\\
&\quad+\sL_1(X_s^\eps,\cL_{X_s^\eps},Y_s^\eps,\cL_{Y_s^\eps})\psi(s,X_s^\eps,\cL_{X_s^\eps},Y_s^\eps,\cL_{Y_s^\eps})\no\\ &\quad+\frac{1}{\eps}\sL_2(X_s^\eps,\cL_{X_s^\eps},Y_s^\eps,\cL_{Y_s^\eps})\psi(s,X_s^\eps,\cL_{X_s^\eps},Y_s^\eps,\cL_{Y_s^\eps})\no\\
&\quad+\frac{1}{\eps}\sL_3(X_s^\eps,\cL_{X_s^\eps},Y_s^\eps,\cL_{Y_s^\eps})\psi(s,X_s^\eps,\cL_{X_s^\eps},Y_s^\eps,\cL_{Y_s^\eps})\no\\ &\quad+\frac{1}{\eps^2}\sL_0(\cL_{X_s^\eps},Y_s^\eps,\cL_{Y_s^\eps}) \psi(s,X_s^\eps,\cL_{X_s^\eps},Y_s^\eps,\cL_{Y_s^\eps})\Big]\dif s+M_t^1+\frac{1}{\eps}M_t^2\no\\
&\quad+\tilde\mE\bigg(\int_0^tF(\tilde X^{\eps}_s,\cL_{X_s^\eps},\tilde Y^{\eps}_s,\cL_{Y^{\eps}_s})\cdot\p_\mu\psi(s,X_s^\eps,\cL_{X_s^\eps},Y_s^\eps,\cL_{Y_s^\eps})(\tilde X^{\eps}_s)\no\\
&\qquad\quad+\frac{1}{\eps}H(\tilde X^{\eps}_s,\cL_{X_s^\eps},\tilde Y^{\eps}_s,\cL_{Y^{\eps}_s})\cdot\p_\mu\psi(s,X_s^\eps,\cL_{X_s^\eps},Y_s^\eps,\cL_{Y_s^\eps})(\tilde X^{\eps}_s)\no\\
&\qquad\quad+\frac{1}{2}\T\Big(GG^*(\tilde X^{\eps}_s,\cL_{X_s^\eps},\tilde Y^{\eps}_s,\cL_{Y^{\eps}_s})\cdot\p_{\tilde x}\big[\p_\mu\psi(s,X_s^\eps,\cL_{X_s^\eps},Y_s^\eps,\cL_{Y_s^\eps})(\tilde X^{\eps}_s)\big]\Big)\no\\
&\qquad\quad+\frac{1}{\eps}c(\tilde X^{\eps}_s,\cL_{X_s^\eps},\tilde Y^{\eps}_s,\cL_{Y^{\eps}_s})\cdot\p_\nu\psi(s,X_s^\eps,\cL_{X_s^\eps},Y_s^\eps,\cL_{Y_s^\eps})(\tilde Y^{\eps}_s)\dif s\bigg),\label{ff1}
\end{align}
where the operators $\sL_1, \sL_2$ and $\sL_3$ are defined by (\ref{op1}), (\ref{op2}) and (\ref{op3}), respectively, the process ($\tilde X^{\eps}_s,\tilde Y^{\eps}_s$) is a copy of the original process $(X^{\eps}_s,Y^{\eps}_s)$ defined on a copy $(\tilde\Omega,\tilde\sF,\tilde\mP)$ of the original probability space $(\Omega,\sF,\mP)$, and $M_t^1$ and $M_t^2$ are two martingales defined by
\begin{align*}
M_t^1&:=\int_0^t\p_x\psi(s,X_s^\eps,\cL_{X_s^\eps},Y_s^\eps,\cL_{Y_s^\eps})\cdot G(X_s^\eps,\cL_{X_s^\eps},Y_s^\eps,\cL_{Y_s^\eps})dW^1_s,\\
M_t^2&:=\int_0^t\p_y\psi(s,X_s^\eps,\cL_{X_s^\eps},Y_s^\eps,\cL_{Y_s^\eps})\cdot\sigma(\cL_{X_s^\eps},Y_s^\eps,\cL_{Y_s^\eps})dW^2_s.
\end{align*}
 Taking expectations  and multiplying $\eps^2$ from both sides of (\ref{ff1}), and in view of  the equation (\ref{pss}), we obtain
\begin{align*}
&\mE\left(\int_0^tf(s,X_s^\eps,\cL_{X_s^\eps},Y_s^\eps,\cL_{Y_s^\eps})\dif s\right)\\
&=\eps^2\,\mE\big[\psi(0,\xi,\cL_\xi,\eta,\cL_\eta)-\psi(t,X_t^\eps,\cL_{X_t^\eps},Y_t^\eps,\cL_{Y_t^\eps})\big]\\
&\quad+\eps^2\,\mE\left(\int_0^t\p_s\psi(s,X_s^\eps,\cL_{X_s^\eps},Y_s^\eps,\cL_{Y_s^\eps})
\dif s\right)\\
&\quad+\eps^2\,\mE\left(\int_0^t
\sL_1(X_s^\eps,\cL_{X_s^\eps},Y_s^\eps,\cL_{Y_s^\eps})\psi(s,X_s^\eps,\cL_{X_s^\eps},Y_s^\eps,\cL_{Y_s^\eps})\dif s\right)\\
&\quad+\eps\,\mE\left(\int_0^t\sL_2(X_s^\eps,\cL_{X_s^\eps},Y_s^\eps,\cL_{Y_s^\eps}) \psi(s,X_s^\eps,\cL_{X_s^\eps},Y_s^\eps,\cL_{Y_s^\eps})\dif s\right)\\
&\quad+\eps\,\mE\left(\int_0^t\sL_3(X_s^\eps,\cL_{X_s^\eps},Y_s^\eps,\cL_{Y_s^\eps}) \psi(s,X_s^\eps,\cL_{X_s^\eps},Y_s^\eps,\cL_{Y_s^\eps})\dif s\right)\\
&\quad+\eps^2\,\mE\tilde\mE\bigg(\int_0^tF(\tilde X^{\eps}_s,\cL_{X_s^\eps},\tilde Y^{\eps}_s,\cL_{Y^{\eps}_s})\cdot\p_\mu\psi(s,X_s^\eps,\cL_{X_s^\eps},Y_s^\eps,\cL_{Y_s^\eps})(\tilde X^{\eps}_s)\no\\
&\qquad\quad+\frac{1}{\eps}H(\tilde X^{\eps}_s,\cL_{X_s^\eps},\tilde Y^{\eps}_s,\cL_{Y^{\eps}_s})\cdot\p_\mu\psi(s,X_s^\eps,\cL_{X_s^\eps},Y_s^\eps,\cL_{Y_s^\eps})(\tilde X^{\eps}_s)\no\\
&\qquad\quad+\frac{1}{2}\T\Big(GG^*(\tilde X^{\eps}_s,\cL_{X_s^\eps},\tilde Y^{\eps}_s,\cL_{Y^{\eps}_s})\cdot\p_{\tilde x}\big[\p_\mu\psi(s,X_s^\eps,\cL_{X_s^\eps},Y_s^\eps,\cL_{Y_s^\eps})(\tilde X^{\eps}_s)\big]\Big)\no\\
&\qquad\quad+\frac{1}{\eps}c(\tilde X^{\eps}_s,\cL_{X_s^\eps},\tilde Y^{\eps}_s,\cL_{Y^{\eps}_s})\cdot\p_\nu\psi(s,X_s^\eps,\cL_{X_s^\eps},Y_s^\eps,\cL_{Y_s^\eps})(\tilde Y^{\eps}_s)\dif s\bigg).
\end{align*}
Using the assumptions on the coefficients and the regularity of $\psi$ again, the expectations on the right hand side of the above equality can be controlled. Thus we arrive at
\begin{align*}
\left|\mE\left(\int_0^tf(s,X_s^\eps,\cL_{X_s^\eps},Y_s^\eps,\cL_{Y_s^\eps})\dif s\right)\right|\leq C_0\,\eps.
\end{align*}
The proof is finished.
\end{proof}

Next, we provide the fluctuations estimate of the functional central limit type. Recall that $\psi$ is the solution of the Poisson equation (\ref{pss}), and $\zeta^\mu$ is the invariant measure for the SDE (\ref{sde1}).
For simplify, we define
\begin{align}
\overline{H\cdot\p_x\psi}(t,x,\mu)&:=\int_{\mR^{d_2}}H(x,\mu,y,\zeta^{\mu})\cdot \p_x\psi(t,x,\mu,y,\zeta^{\mu})\zeta^{\mu}(\dif y),\label{hpsi}\\
\overline{c\cdot\p_y\psi}(t,x,\mu)&:=\int_{\mR^{d_2}}c(x,\mu,y,\zeta^{\mu})\cdot \p_y\psi(t,x,\mu,y,\zeta^{\mu})\zeta^{\mu}(\dif y),\label{cpsi1}
\end{align}
and
\begin{align}
\overline{c\cdot\p_\nu\psi}(t,x,\mu,y)(\tilde x)&:=\int_{\mR^{d_2}}c(\tilde x,\mu,\tilde y,\zeta^{\mu})\cdot \p_\nu\psi(t,x,\mu,y,\zeta^{\mu})(\tilde y)\zeta^{\mu}(\dif \tilde y),\label{cpsi2}\\
\overline{\overline{c\cdot\p_\nu\psi}}(t,x,\mu)(\tilde x)&:=\int_{\mR^{d_2}} \overline{c\cdot\p_\nu\psi}(t,x,\mu,y)(\tilde x)\zeta^{\mu}(\dif y)\no\\
&=\int_{\mR^{d_2}}\int_{\mR^{d_2}}c(\tilde x,\mu,\tilde y,\zeta^{\mu})\cdot \p_\nu\psi(t,x,\mu,y,\zeta^{\mu})(\tilde y)\zeta^{\mu}(\dif \tilde y)\zeta^{\mu}(\dif y).\label{cpsi3}
\end{align}
The following result will play an important role below.

\bl\label{flu2}
Assume that {\bf (H$^{\sigma,b}$)} hold, $F,H,G,c\in C_b^{2,(1,1),2,(1,1)}(\mR^{d_1}\times\sP_2(\mR^{d_1})\times\mR^{d_2}\times\sP_2(\mR^{d_2}))$ and
$\sigma, b\in \mC_b^{(1,2),2,(1,1)}\cap C_b^{(2,1),4,(2,2)}(\sP_2(\mR^{d_1})\times\mR^{d_2}\times\sP_2(\mR^{d_2}))$. Then for every $f\in \mC_b^{1,3,(1,2),2,(1,1)}\cap C_b^{1,3,(2,1),4,(2,2)}(\mR_+\times\mR^{d_1}\times\sP_2(\mR^{d_1})\times\mR^{d_2}\times\sP_2(\mR^{d_2}))$
satisfying (\ref{cen}), we have
\begin{align}\label{es2}
&\bigg|\mE\left(\frac{1}{\eps}\int_0^tf(s,X_s^\eps,\cL_{X_s^\eps},Y_s^\eps,\cL_{Y_s^\eps})\dif s\right)\no\\
&\qquad-\mE\bigg(\int_0^t\overline{H\cdot\p_x\psi}(s,X_s^\eps,\cL_{X_s^\eps})
+\overline{c\cdot\p_y\psi}(s,X_s^\eps,\cL_{X_s^\eps})\dif s\bigg)\no\\
&\qquad-\mE\tilde{\mE}\bigg(\int_0^t\overline{\overline{c\cdot\p_\nu\psi}}(s,X_s^\eps,\cL_{X_s^\eps})(\tilde{X}_s^\eps)\dif s\bigg)\bigg|\leq C_0\,\eps,
\end{align}
where $C_0>0$ is a constant independent of $\eps$.
\el
\begin{proof}
Multiplying $\eps$ from both sides of (\ref{ff1}) and using (\ref{pss}), we have
\begin{align*}
&\mE\left(\frac{1}{\eps}\int_0^tf(s,X_s^\eps,\cL_{X_s^\eps},Y_s^\eps,\cL_{Y_s^\eps})\dif s\right)\\
&=\eps\,\mE\big[\psi(0,\xi,\cL_\xi,\eta,\cL_\eta)-\psi(t,X_t^\eps,\cL_{X_t^\eps},Y_t^\eps,\cL_{Y_t^\eps})\big]\\
&\quad+\eps\,\mE\left(\int_0^t\p_s\psi(s,X_s^\eps,\cL_{X_s^\eps},Y_s^\eps,\cL_{Y_s^\eps})\dif s\right)\\
&\quad+\eps\,\mE\left(\int_0^t
\sL_1(X_s^\eps,\cL_{X_s^\eps},Y_s^\eps,\cL_{Y_s^\eps})\psi(s,X_s^\eps,\cL_{X_s^\eps},Y_s^\eps,\cL_{Y_s^\eps})\dif s\right)\\
&\quad+\eps\,\mE\tilde\mE\bigg(\int_0^tF(\tilde X^{\eps}_s,\cL_{X_s^\eps},\tilde Y^{\eps}_s,\cL_{Y^{\eps}_s})\cdot\p_\mu\psi(s,X_s^\eps,\cL_{X_s^\eps},Y_s^\eps,\cL_{Y_s^\eps})(\tilde X^{\eps}_s)\no\\
&\qquad\quad+\frac{1}{2}\T\Big(GG^*(\tilde X^{\eps}_s,\cL_{X_s^\eps},\tilde Y^{\eps}_s,\cL_{Y^{\eps}_s})\cdot\p_{\tilde x}\big[\p_\mu\psi(s,X_s^\eps,\cL_{X_s^\eps},Y_s^\eps,\cL_{Y_s^\eps})(\tilde X^{\eps}_s)\big]\Big)\dif s\bigg)\no\\
&\quad+\mE\left(\int_0^t H(X_s^\eps,\cL_{X_s^\eps},Y_s^\eps,\cL_{Y_s^\eps})\cdot \p_x\psi(s,X_s^\eps,\cL_{X_s^\eps},Y_s^\eps,\cL_{Y_s^\eps})\dif s\right)\\
&\quad+\mE\left(\int_0^t c(X_s^\eps,\cL_{X_s^\eps},Y_s^\eps,\cL_{Y_s^\eps})\cdot \p_y\psi(s,X_s^\eps,\cL_{X_s^\eps},Y_s^\eps,\cL_{Y_s^\eps})\dif s\right)\\
&\quad+\mE\tilde\mE\left(\int_0^t H(\tilde X^{\eps}_s,\cL_{X_s^\eps},\tilde Y^{\eps}_s,\cL_{Y^{\eps}_s})\cdot\p_\mu\psi(s,X_s^\eps,\cL_{X_s^\eps},Y_s^\eps,\cL_{Y_s^\eps})(\tilde X^{\eps}_s)\dif s\right)\no\\
&\quad+\mE\tilde\mE\left(\int_0^t c(\tilde X^{\eps}_s,\cL_{X_s^\eps},\tilde Y^{\eps}_s,\cL_{Y^{\eps}_s})\cdot\p_\nu\psi(s,X_s^\eps,\cL_{X_s^\eps},Y_s^\eps,\cL_{Y_s^\eps})(\tilde Y^{\eps}_s)\dif s\right).
\end{align*}
By the same argument as in the proof of Lemma \ref{flu1}, we obtain
\begin{align*}
&\bigg|\mE\left(\frac{1}{\eps}\int_0^tf(s,X_s^\eps,\cL_{X_s^\eps},Y_s^\eps,\cL_{Y_s^\eps})\dif s\right)-\mE\bigg(\int_0^t\overline{H\cdot\p_x\psi}(s,X_s^\eps,\cL_{X_s^\eps})\dif s\bigg)\\
&-\mE\bigg(\int_0^t\overline{c\cdot\p_y\psi}(s,X_s^\eps,\cL_{X_s^\eps})\dif s\bigg)
-\mE\tilde{\mE}\bigg(\int_0^t\overline{\overline{c\cdot\p_\nu\psi}}(s,X_s^\eps,\cL_{X_s^\eps})(\tilde{X}_s^\eps)\dif s\bigg)\bigg|\\
&\leq C_0\,\eps+\bigg|\mE\tilde\mE\left(\int_0^tH(\tilde X^{\eps}_s,\cL_{X_s^\eps},\tilde Y^{\eps}_s,\cL_{Y^{\eps}_s})\cdot\p_\mu\psi(s,X_s^\eps,\cL_{X_s^\eps},Y_s^\eps,\cL_{Y_s^\eps})(\tilde X^{\eps}_s)\dif s\right)\bigg|\\
&\quad+\bigg|\mE\left(\int_0^t H(X_s^\eps,\cL_{X_s^\eps},Y_s^\eps,\cL_{Y_s^\eps})\cdot \p_x\psi(s,X_s^\eps,\cL_{X_s^\eps},Y_s^\eps,\cL_{Y_s^\eps})\dif s\right)\\
&\qquad-\mE\bigg(\int_0^t\overline{H\cdot\p_x\psi}(s,X_s^\eps,\cL_{X_s^\eps})\dif s\bigg)\bigg|\\
&\quad+\bigg|\mE\left(\int_0^t c(X_s^\eps,\cL_{X_s^\eps},Y_s^\eps,\cL_{Y_s^\eps})\cdot \p_y\psi(s,X_s^\eps,\cL_{X_s^\eps},Y_s^\eps,\cL_{Y_s^\eps})\dif s\right)\\
&\qquad-\mE\bigg(\int_0^t\overline{c\cdot\p_y\psi}(s,X_s^\eps,\cL_{X_s^\eps})\dif s\bigg)\bigg|\\
&\quad+\bigg|\mE\tilde \mE\left(\int_0^t c(\tilde X_s^\eps,\cL_{X_s^\eps},\tilde Y_s^\eps,\cL_{Y_s^\eps})\cdot \p_\nu\psi(s,X_s^\eps,\cL_{X_s^\eps},Y_s^\eps,\cL_{Y_s^\eps})(\tilde Y_s^\eps)\dif s\right)\\
&\qquad-\mE\tilde \mE\bigg(\int_0^t\overline{\overline{c\cdot\p_\nu\psi}}(s,X_s^\eps,\cL_{X_s^\eps})(\tilde X_s^\eps)\dif s\bigg)\bigg|\\
&=:C_0\,\eps+\cI_1(\eps)+\cI_2(\eps)+\cI_3(\eps)+\cI_4(\eps).
\end{align*}
In what follows, we estimate the above four terms one by one. For the first term, we write
\begin{align*}
\cI_1(\eps)=\bigg|\mE\bigg[\tilde\mE\left(\int_0^tH(\tilde X^{\eps}_s,\cL_{X_s^\eps},\tilde Y^{\eps}_s,\cL_{Y^{\eps}_s})\cdot\p_\mu\psi(s,x,\cL_{X_s^\eps},y,\cL_{Y_s^\eps})(\tilde X^{\eps}_s)\dif s\right)\bigg|_{(x,y)=(X^{\eps}_s,Y^{\eps}_s)}\bigg]\bigg|.
\end{align*}
Since $H(\tilde x,\mu,\tilde y,\nu)$ satisfies the centering condition (\ref{cen}), this in turn implies that for every fixed $(x,y)\in\mR^{d_1}\times\mR^{d_2}$,
$$
(t,\tilde x,\mu,\tilde y,\nu)\mapsto H(\tilde x,\mu,\tilde y,\nu)\p_\mu\psi(t,x,\mu,y,\nu)(\tilde x)
$$
satisfies the centering condition, too. Moreover, by the assumptions on the coefficients and Theorem \ref{po}, we have
$$
\p_\mu\psi(\cdot,x,\cdot,y,\cdot)(\cdot)\in C_b^{1,(1,1),(1,1),2}(\mR_+\times\sP_2(\mR^{d_1})\times\sP_2(\mR^{d_2})\times\mR^{d_1}).
$$
Thus using Lemma \ref{flu1} we get that for every fixed $(x,y)\in\mR^{d_1}\times\mR^{d_2}$,
\begin{align*}
\bigg|\tilde\mE\left(\int_0^tH(\tilde X^{\eps}_s,\cL_{X_s^\eps},\tilde Y^{\eps}_s,\cL_{Y^{\eps}_s})\cdot\p_\mu\psi(s,x,\cL_{X_s^\eps},y,\cL_{Y_s^\eps})(\tilde X^{\eps}_s)\dif s\right)\bigg|\leq C_1 \,\eps,
\end{align*}
which in turn implies that
$$
\cI_1(\eps)\leq C_1\,\eps.
$$
To control the second and third term, note that by the definition  (\ref{hpsi}) and (\ref{cpsi1}) we have
\begin{align*}
\int_{\mathbb{R}^{d_2}} \big[H(x,\mu,y,\zeta^\mu)\cdot\p_x\psi(t,x,\mu,y,\zeta^\mu)-\overline{H\cdot\p_x\psi}(t,x,\mu)\big]\zeta^\mu(\dif y)=0
\end{align*}
and
\begin{align*}
\int_{\mathbb{R}^{d_2}} \big[c(x,\mu,y,\zeta^\mu)\cdot\p_y\psi(t,x,\mu,y,\zeta^\mu)-\overline{c\cdot\p_y\psi}(t,x,\mu)\big]\zeta^\mu(\dif y)=0.
\end{align*}
Since $\p_x\psi, \p_y\psi \in C_b^{1,2,(1,1),2,(1,1)}(\mR_+\times\mR^{d_1}\times\sP_2(\mR^{d_1})\times\mR^{d_2}\times\sP_2(\mR^{d_2}))$, we have by Corollary \ref{avef} that
$$
\overline{H\cdot\p_x\psi}(t,x,\mu), \overline{c\cdot\p_y\psi}(t,x,\mu)\in C_b^{1,2,(1,1)}(\mR_+\times\mR^{d_1}\times\sP_2(\mR^{d_1})).
$$
Using Lemma \ref{flu1} again we obtain
$$
\cI_2(\eps)+\cI_3(\eps)\leq C_2\,\eps.
$$
As for $\cI_4(\eps)$, using (\ref{cpsi2}) we write
\begin{align}
\cI_4(\eps)&\leq\bigg|\mE\tilde \mE\left(\int_0^t c(\tilde X_s^\eps,\cL_{X_s^\eps},\tilde Y_s^\eps,\cL_{Y_s^\eps})\cdot \p_\nu\psi(s,X_s^\eps,\cL_{X_s^\eps},Y_s^\eps,\cL_{Y_s^\eps})(\tilde Y_s^\eps)\dif s\right)\no\\
&\qquad-\mE\tilde \mE\left(\int_0^t \overline{c\cdot\p_\nu\psi}(s,X_s^\eps,\cL_{X_s^\eps},Y_s^\eps)(\tilde X_s^\eps)\dif s\right)\bigg|\no\\
&\quad+\bigg|\mE\tilde \mE\left(\int_0^t \overline{c\cdot\p_\nu\psi}(s,X_s^\eps,\cL_{X_s^\eps},Y_s^\eps)(\tilde X_s^\eps)\dif s\right)\no\\
&\qquad-\mE\tilde \mE\bigg(\int_0^t\overline{\overline{c\cdot\p_\nu\psi}}(s,X_s^\eps,\cL_{X_s^\eps})(\tilde X_s^\eps)\dif s\bigg)\bigg|\no\\
&\leq\bigg|\mE\bigg[\tilde \mE\bigg(\int_0^t c(\tilde X_s^\eps,\cL_{X_s^\eps},\tilde Y_s^\eps,\cL_{Y_s^\eps})\cdot \p_\nu\psi(s,x,\cL_{X_s^\eps},y,\cL_{Y_s^\eps})(\tilde Y_s^\eps)\no\\
&\qquad-\overline{c\cdot\p_\nu\psi}(s,x,\cL_{X_s^\eps},y)(\tilde X_s^\eps)\dif s\bigg)\bigg|_{(x,y)=(X_s^\eps,Y_s^\eps)}\bigg]\bigg|\no\\
&\quad+\bigg|\tilde \mE\bigg[\mE\bigg(\int_0^t \overline{c\cdot\p_\nu\psi}(s,X_s^\eps,\cL_{X_s^\eps},Y_s^\eps)(\tilde x)
-\overline{\overline{c\cdot\p_\nu\psi}}(s,X_s^\eps,\cL_{X_s^\eps})(\tilde x)\dif s\bigg)\bigg|_{\tilde x=\tilde X_s^\eps}\bigg]\bigg|\no\\
&=:\cI_{4,1}(\eps)+\cI_{4,2}(\eps).\label{estc4}
\end{align}
By the definition of $\overline{c\cdot\p_\nu\psi}(t,x,\mu,y)(\tilde x)$, for any fixed $(x,y)\in\mathbb{R}^{d_1}\times\mathbb{R}^{d_2}$,
$$
(\tilde x,\mu,\tilde y,\nu)\mapsto c(\tilde x,\mu,\tilde y,\nu)\cdot\p_\nu\psi(t,x,\mu,y,\nu)(\tilde y)-\overline{c\cdot\p_\nu\psi}(t,x,\mu,y)(\tilde x)
$$
satisfies the centering condition (\ref{cen}). By the assumptions on the coefficients and Theorem \ref{po}, we have
$$
\p_\nu\psi\in C_b^{1,2,(1,1),2,(1,1),2}(\mR_+\times\mR^{d_1}\times\sP_2(\mR^{d_1})\times\mR^{d_2}\times\sP_2(\mR^{d_2})\times\mR^{d_2}).
$$
 Therefore, using Corollary \ref{avef} and  Lemma \ref{flu1} we obtain
\begin{align*}
&\bigg|\tilde \mE\bigg(\int_0^t c(\tilde X_s^\eps,\cL_{X_s^\eps},\tilde Y_s^\eps,\cL_{Y_s^\eps})\cdot \p_\nu\psi(s,x,\cL_{X_s^\eps},y,\cL_{Y_s^\eps})(\tilde Y_s^\eps)\no\\
&\quad-\overline{c\cdot\p_\nu\psi}(s,x,\cL_{X_s^\eps},y)(\tilde X_s^\eps)\dif s\bigg)\bigg|\leq C_4 \,\eps,
\end{align*}
which in turn implies
\begin{align}\label{estc41}
\cI_{4,1}(\eps)\leq C_4\,\eps.
\end{align}
Similarly, in view of (\ref{cpsi3}) we have that for any fixed $\tilde x\in\mathbb{R}^{d_1}$,
$$
(x,\mu,y)\mapsto \overline{c\cdot\p_\nu\psi}(t,x,\mu,y)(\tilde x)-\overline{\overline{c\cdot\p_\nu\psi}}(t,x,\mu)(\tilde x)
$$
also satisfies the centering condition. By the same argument as (\ref{estc41}) we have
\begin{align}\label{estc42}
\cI_{4,2}(\eps)\leq C_4\,\eps.
\end{align}
Substituting (\ref{estc41}) and (\ref{estc42}) into (\ref{estc4}) yields
$$
\cI_{4}(\eps)\leq C_4\,\eps.
$$
Combing the above computations, the proof is finished.
\end{proof}

\section{Proof of Theorem \ref{main}}

Throughout this section, we assume that the conditions in Theorem \ref{main} hold.
Let $\bar X^{s,\xi}_t$ be the unique solution to the SDE (\ref{ave}) starting from the initial point $\xi\in L^2(\Omega)$ at time $s$.
Namely, for $t\geq s$,
\begin{align*}
\dif \bar X^{s,\xi}_t&=\bar F(\bar X^{s,\xi}_t,\cL_{\bar X^{s,\xi}_t})\dif t+\overline{H\cdot\p_x\Phi}(\bar X^{s,\xi}_t,\cL_{\bar X^{s,\xi}_t})\dif t\\
&\quad+\overline{c\cdot\p_y\Phi}(\bar X^{s,\xi}_t,\cL_{\bar X^{s,\xi}_t})\dif t+\int_{\mathbb{R}^{d_1}}\overline{\overline{c\cdot\p_\nu\Phi}}(\bar X^{s,\xi}_t,\cL_{\bar X^{s,\xi}_t})(\tilde x)\cL_{\bar X^{s,\xi}_t}(\dif \tilde x)\dif t\no\\
&\quad+\sqrt{\overline{GG^*}+2\overline{H\cdot\Phi}}(\bar X^{s,\xi}_t,\cL_{\bar X^{s,\xi}_t})\dif W^1_t,\quad  \bar X_s^{s,\xi}=\xi,
\end{align*}
where the coefficients are defined by (\ref{bF})-(\ref{bc3}), respectively. Fix $T>0$ and $\varphi: \sP_2(\mR^{d_1})\to\mR$. For $t\in[0,T]$, define
\begin{align}\label{uu}
u(t,\cL_\xi):=\varphi(\cL_{\bar X^{t,\xi}_T}).
\end{align}
Then we have:
\bl
Assume that $\varphi\in C_b^{(3,1)}(\sP_2(\mR^{d_1}))$. Then $u(t,\cL_\xi)$ is the unique solution in $C_b^{1,(3,1)}([0,T]\times\sP_2(\mR^{d_1}))$ of the equation
\begin{equation}\label{equ}
\left\{ \begin{aligned}
&\p_t u(t,\cL_\xi)+\mE\Big[\Big(\bar F(\xi,\cL_{\xi})+\overline{H\cdot\p_x\Phi}(\xi,\cL_{\xi}) +\overline{c\cdot\p_y\Phi}(\xi,\cL_{\xi})\Big)\cdot\p_\mu u(t,\cL_{\xi})(\xi)\\
&\qquad\qquad+\int_{\mathbb{R}^{d_1}}\overline{\overline{c\cdot\p_\nu\Phi}}(\xi,\cL_{\xi})(\tilde x)\cL_{\xi}(\dif \tilde x)\cdot\p_\mu u(t,\cL_{\xi})(\xi)\\
&\qquad\qquad+\frac{1}{2}\T\Big(\big(\overline{GG^*}+2\overline{H\cdot\Phi}(\xi,\cL_{\xi})\big)\cdot\p_{x}\big[\p_\mu u(t,\cL_{\xi})\big](\xi)\Big)\Big]=0,\\
&u(T,\cL_{\xi})=\varphi(\cL_\xi).
\end{aligned} \right.
\end{equation}
\el
\begin{proof}
Under the assumptions on the coefficients and using Corollary \ref{avef}, we have
$\bar F,\overline{GG^*}\in C_b^{3,(3,1)}(\mR^{d_1}\times\sP_2(\mR^{d_1}))$. Meanwhile, since $\Phi\in \textbf{C}_b^{4,(2,2),4,(2,2)}\cap \cC_b^{4,6,(3,3)}(\mR^{d_1}\times\sP_2(\mR^{d_1})\times\mR^{d_2}\times\sP_2(\mR^{d_2}))$, using Corollary \ref{avef} again we have $\overline{H\cdot\p_x\Phi}, \overline{c\cdot\p_y\Phi}, \overline{H\cdot\Phi}\in C_b^{3,(3,1)}(\mR^{d_1}\times\sP_2(\mR^{d_1}))$ and $\overline{\overline{c\cdot\p_\nu\Phi}}(x,\mu)(\tilde x)\in C_b^{3,(3,1),2}(\mR^{d_1}\times\sP_2(\mR^{d_1})\times\mR^{d_1})$.
The statements follows by \cite[Theorem 7.2]{BLPR} and the same arguments as in the proof of Theorem \ref{pot}, we omit the details.
\end{proof}

We are now in the position to give:

\begin{proof}[Proof of Theorem \ref{main}]
Let $u(t,\cL_\xi)$ be defined by (\ref{uu}). Then we have
\begin{align*}
\sJ(\eps):=\varphi(\cL_{X_T^\eps})-\varphi(\cL_{\bar X_T})=u(T,\cL_{X_T^\eps})-u(0,\cL_\xi).
\end{align*}
Thus, by It\^o's formula,
\begin{align*}
\sJ(\eps)&=\mE\bigg(\int_0^T\p_tu(t,\cL_{X_t^\eps}) +F(X^{\eps}_t,\cL_{X_t^\eps},Y^{\eps}_t,\cL_{Y^{\eps}_t})\cdot\p_{\mu}u(t,\cL_{X_t^\eps})(X^{\eps}_t)\\
&\quad+\frac{1}{\eps}H(X^{\eps}_t,\cL_{X_t^\eps},Y^{\eps}_t,\cL_{Y^{\eps}_t}) \cdot\p_{\mu}u(t,\cL_{X_t^\eps})(X^{\eps}_t)\\ &\quad+\frac{1}{2}\T\Big(GG^*(X^{\eps}_t,\cL_{X_t^\eps},Y^{\eps}_t,\cL_{Y^{\eps}_t})
\cdot\p_{x}\big[\p_{\mu}u(t,\cL_{X_t^\eps})\big](X^{\eps}_t)\Big)\dif t\bigg).
\end{align*}
In view of equation (\ref{equ}), we further obtain that
\begin{align*}
\big|\sJ(\eps)\big|&\leq\bigg|\mE\bigg(\int_0^T\delta F(X^{\eps}_t,\cL_{X_t^\eps},Y^{\eps}_t,\cL_{Y^{\eps}_t})\cdot\p_{\mu}u(t,\cL_{X_t^\eps})(X^{\eps}_t)\dif t\bigg)\bigg|\\
&\quad+\frac{1}{2}\bigg|\mE\bigg(\int_0^T\T\Big(\delta (GG^*)(X^{\eps}_t,\cL_{X_t^\eps},Y^{\eps}_t,\cL_{Y^{\eps}_t})\cdot\p_{ x}\big[\p_{\mu}u(t,\cL_{X_t^\eps})\big](X^{\eps}_t)\Big)\dif t\bigg)\bigg|\\
&\quad+\bigg|\mE\bigg(\int_0^T\frac{1}{\eps}H(X^{\eps}_t,\cL_{X_t^\eps},Y^{\eps}_t,\cL_{Y^{\eps}_t}) \cdot\p_{\mu}u(t,\cL_{X_t^\eps})(X^{\eps}_t)\\
&\quad-\Big(\overline{H\cdot\p_x\Phi}(X_t^\eps,\cL_{X_t^\eps}) +\overline{c\cdot\p_y\Phi}(X_t^\eps,\cL_{X_t^\eps})\Big)\cdot\p_\mu u(t,\cL_{X_t^\eps})(X_t^\eps)\\
&\quad-\int_{\mathbb{R}^{d_1}}\overline{\overline{c\cdot\p_\nu\Phi}}(X_t^\eps,\cL_{X_t^\eps})(\tilde x)\cL_{X_t^\eps}(\dif \tilde x)\cdot\p_\mu u(t,\cL_{X_t^\eps})(X_t^\eps)\\
&\quad-\T\Big(\overline{H\cdot\Phi}(X_t^\eps,\cL_{X_t^\eps})\cdot\p_{x}\big[\p_\mu u(t,\cL_{X_t^\eps})\big](X_t^\eps)\Big)\dif t\bigg)\bigg|\\
&=:\sJ_1(\eps)+\sJ_2(\eps)+\sJ_3(\eps),
\end{align*}
where
$$
\delta F(x,\mu,y,\nu):=F(x,\mu,y,\nu)-\bar F(x,\mu)
$$
and
$$
\delta(GG^*)(x,\mu,y,\nu):=GG^*(x,\mu,y,\nu)-\overline{GG^*}(x,\mu).
$$
By the definition of $\bar F(x,\mu)$, we have
\begin{align*}
&\int_{\mR^{d_2}}\delta F(x,\mu,y,\zeta^\mu)\cdot\p_{\mu}u(t,\mu)(x)\zeta^{\mu}(\dif y)\\
&=\int_{\mR^{d_2}}\delta F(x,\mu,y,\zeta^\mu)\zeta^{\mu}(\dif y)\cdot\p_{\mu}u(t,\mu)(x)=0.
\end{align*}
Meanwhile, $\delta F(x,\mu,y,\nu)\cdot\p_{\mu}u(t,\mu)(x)\in C_b^{1,2,(1,1),2,(1,1)}(\mR_+\times\mR^{d_1}\times\sP_2(\mR^{d_1})\times\mR^{d_2}\times\sP_2(\mR^{d_2}))$. As a result of Lemma \ref{flu1}, we have
$$
\sJ_1(\eps)\leq C_1\,\eps.
$$
Note that
\begin{align*}
&\int_{\mR^{d_2}}\delta(GG^*)(x,\mu,y,\zeta^\mu)\cdot\p_{x}[\p_{\mu}u(t,\mu)(x)]\zeta^{\mu}(\dif y)\\
&=\int_{\mR^{d_2}}\delta(GG^*)(x,\mu,y,\zeta^\mu)\zeta^{\mu}(\dif y)\cdot\p_{x}[\p_{\mu}u(t,\mu)(x)]=0.
\end{align*}
Then, using Lemma \ref{flu1} again, one can get
$$
\sJ_2(\eps)\leq C_2\,\eps.
$$
To control the third term, note that since $H(x,\mu,y,\nu)$ satisfies the centering condition (\ref{cen}), we have
\begin{align*}
&\int_{\mR^{d_2}}H(x,\mu,y,\zeta^\mu)\cdot\p_{\mu}u(t,\mu)(x)\zeta^{\mu}(\dif y)\\
&=\int_{\mR^{d_2}}H(x,\mu,y,\zeta^\mu)\zeta^{\mu}(\dif y)\cdot\p_{\mu}u(t,\mu)(x)=0.
\end{align*}
On the other hand, recall that $\Phi(x,\mu,y,\nu)$ satisfies the Poisson equation (\ref{po1}), and define
$$
\tilde\Phi(t,x,\mu,y,\nu):=\Phi(x,\mu,y,\nu)\cdot\p_{\mu}u(t,\mu)(x).
$$
Then, one can check that
$$
\sL_0\tilde\Phi(t,x,\mu,y,\nu)=-H(x,\mu,y,\nu)\cdot\p_{\mu}u(t,\mu)(x).
$$
Moreover, we have $\tilde\Phi\in C_b^{1,2,(1,1),2,(1,1)}(\mR_+\times\mR^{d_1}\times\sP_2(\mR^{d_1})\times\mR^{d_2}\times\sP_2(\mR^{d_2}))$ and it holds that
\begin{align*}
\overline{H\cdot\p_x\tilde\Phi}(t,x,\mu)&=\int_{\mR^{d_2}} H(x,\mu,y,\zeta^{\mu})\cdot\p_x\Phi(x,\mu,y,\zeta^{\mu})\cdot\p_{\mu}u(t,\mu)(x)\zeta^{\mu}(\dif y)\\
&\quad+\int_{\mR^{d_2}}\T \Big(H(x,\mu,y,\zeta^{\mu})\cdot\Phi(x,\mu,y,\zeta^{\mu})\cdot\p_x\big[\p_{\mu}u(t,\mu)\big](x)\Big)\zeta^{\mu}(\dif y)\\
&=\overline{H\cdot\p_x\Phi}(x,\mu)\cdot\p_{\mu}u(t,\mu)(x) +\T\Big(\overline{H\cdot\Phi}(x,\mu)\cdot\p_x\big[\p_{\mu}u(t,\mu)\big](x)\Big),\\
\overline{c\cdot\p_y\tilde\Phi}(t,x,\mu)&=\int_{\mR^{d_2}} c(x,\mu,y,\zeta^{\mu})\cdot\p_y\Phi(x,\mu,y,\zeta^{\mu})\cdot\p_{\mu}u(t,\mu)(x)\zeta^{\mu}(\dif y)\\
&=\overline{c\cdot\p_y\Phi}(x,\mu)\cdot\p_{\mu}u(t,\mu)(x),
\end{align*}
and
\begin{align*}
&\overline{\overline{H\cdot\p_\nu\tilde\Phi}}(t,x,\mu)(\tilde x)\\
&=\int_{\mR^{d_2}}\int_{\mR^{d_2}} c(\tilde x,\mu,\tilde y,\zeta^{\mu})\cdot\p_\nu\Phi(x,\mu,y,\zeta^{\mu})(\tilde y)\cdot\p_{\mu}u(t,\mu)(x)\zeta^{\mu}(\dif \tilde y)\zeta^{\mu}(\dif y)\\
&=\overline{\overline{c\cdot\p_\nu\Phi}}(x,\mu)(\tilde x)\cdot\p_{\mu}u(t,\mu)(x).
\end{align*}
Thus, by estimate (\ref{es2}) we have
$$
\sJ_3(\eps)\leq C_3\,\eps.
$$
The proof is completed.
\end{proof}

\bigskip

\section{Appendix}

Recall that $Y_t^\eta$ and $Y_t^{y,\nu}$ satisfy the equations (\ref{eqy}) and (\ref{eqyy}), respectively. Throughout this section, we assume the assumption {\bf ($\hat H^{\sigma,b}$)} holds, and provide the following estimates for $Y_t^{y,\nu}$.

\bl\label{Y1}
Assume that {\bf ($\hat H^{\sigma,b}$)} holds. Then for any $p\geq2$, we have
\begin{align*}
&\mE\|\p_yY^{y,\nu}_t\|^p\leq C_0\,\e^{-\frac{p}{2}c_2t},\\
&\mE\|\p^2_yY^{y,\nu}_t\|^p\leq C_0\,\e^{-\frac{p}{2}(c_2-\gamma)t},
\end{align*}
where  $C_0$ is a positive constant independent of $t$ and $\gamma\ll(c_2-c_1)$.
\el
\begin{proof}
     Recall that
	\begin{align*}
	\dif \p_y Y^{y,\nu}_t=\p_y b(Y^{y,\nu}_t,\cL_{Y^\eta_t})\cdot\p_y Y^{y,\nu}_t\dif t+\p_y \sigma(Y^{y,\nu}_t,\cL_{Y^\eta_t})\cdot\p_y Y^{y,\nu}_t \dif W_t.
	\end{align*}
	Using  It\^{o}'s formula, we compute that
	\begin{align}\label{Yy}
	\dif \mE\|\p_y Y^{y,\nu}_t\|^p
    &\leq\frac{p}{2}\mE\big[\|\p_y Y^{y,\nu}_t\|^{p-2}\cdot\big(2\langle \p_y Y^{y,\nu}_t,\p_y b(Y^{y,\nu}_t,\cL_{Y^{\eta}_t})\cdot\p_y Y^{y,\nu}_t \rangle\no\\
	&\quad+(p-1)\|\p_y\sigma(Y^{y,\nu}_t,\cL_{Y^{\eta}_t})\cdot\p_y Y^{y,\nu}_t\|^2\big)\big]\dif t.
	\end{align}
	In view of the assumption {\bf ($\hat H^{\sigma,b}$)}, we have for any $h\in\mR^{d_2}$,
	\begin{align*}
	2\langle h,\p_yb(y,\nu)\cdot h\rangle+(p-1)\|\p_y\sigma(y,\nu)\cdot h\|^2\leq-c_2|h|^2,
	\end{align*}
	 which together with (\ref{Yy}) implies that
	\begin{align*}
	\dif \mE\|\p_y Y^{y,\nu}_t\|^p
	\leq-\frac{p}{2}c_2\mE\|\p_y Y^{y,\nu}_t\|^p\dif t.
	\end{align*}
	Thus, by the comparison theorem, we get
	\begin{align*}
	\mE\|\p_y Y^{y,\nu}_t\|^p\leq C_0e^{-\frac{p}{2}c_2t}.
	\end{align*}
    Similarly, we have
    \begin{align*}
	\dif \mE\|\p^2_y Y^{y,\nu}_t\|^p
    &\leq\frac{p}{2}\mE\big[\|\p^2_y Y^{y,\nu}_t\|^{p-2}\cdot2\langle \p^2_y Y^{y,\nu}_t,\p_y b(Y^{y,\nu}_t,\cL_{Y^{\eta}_t})\cdot\p^2_y Y^{y,\nu}_t\\
    &\qquad+\p^2_y b(Y^{y,\nu}_t,\cL_{Y^{\eta}_t})\cdot\p_y Y^{y,\nu}_t\cdot\p_y Y^{y,\nu}_t \rangle\big]\dif t\no\\
	&\quad+\frac{p(p-1)}{2}\mE\big[\|\p^2_y Y^{y,\nu}_t\|^{p-2}\cdot\|\p_y\sigma(Y^{y,\nu}_t,\cL_{Y^{\eta}_t})\cdot\p^2_y Y^{y,\nu}_t\\
    &\qquad+\p^2_y \sigma(Y^{y,\nu}_t,\cL_{Y^{\eta}_t})\cdot\p_y Y^{y,\nu}_t\cdot\p_y Y^{y,\nu}_t\|^2\big]\dif t\\
    &\leq-\frac{p}{2}(c_2-\gamma)\mE\|\p^2_y Y^{y,\nu}_t\|^p\dif t+C_0\mE\|\p_y Y^{y,\nu}_t\|^{2p}\dif t\\
    &\leq-\frac{p}{2}(c_2-\gamma)\mE\|\p^2_y Y^{y,\nu}_t\|^p\dif t+C_0e^{-pc_2t}\dif t,
	\end{align*}
    which in turn yields \begin{align*}
	\mE\|\p^2_y Y^{y,\nu}_t\|^p\leq C_0e^{-\frac{p}{2}(c_2-\gamma)t}.
	\end{align*}
    Thus the proof is completed.
\end{proof}

\bl\label{Y2}
Assume that {\bf ($\hat H^{\sigma,b}$)} holds. Then we have for every $p\geq2$,
\begin{align*}
&\mE\|\p_\nu Y^{y,\nu}_t(\tilde y)\|^p\leq C_0\,\e^{-\frac{p}{2}(c_2-c_1-\gamma)t},\\
&\mE\|\p_{\tilde y}\p_\nu Y^{y,\nu}_t(\tilde y)\|^p\leq C_0\,\e^{-\frac{p}{2}(c_2-c_1-\gamma)t},\\
&\mE\|\p_y\p_\nu Y^{y,\nu}_t(\tilde y)\|^p\leq C_0\,\e^{-\frac{p}{2}(c_2-c_1-\gamma)t},\\
&\mE\|\p^2_\nu Y^{y,\nu}_t(\tilde y,\breve{y})\|^p\leq C_0\,\e^{-\frac{p}{2}(c_2-c_1-\gamma)t},
\end{align*}
where  $C_0$ is a positive constant independent of $t$ and $\gamma\ll(c_2-c_1)$.
\el
\begin{proof}
Recall that
	\begin{align*}
    \dif \p_\nu Y^{y,\nu}_t(\tilde y)&=\p_y b(Y^{y,\nu}_t,\cL_{Y^\eta_t})\cdot\p_\nu Y^{y,\nu}_t(\tilde y)\dif t+\tilde\mE\big[\p_\nu b(Y^{y,\nu}_t,\cL_{Y^\eta_t})(\tilde Y^{\tilde y,\nu}_t)\cdot\p_y \tilde Y^{\tilde y,\nu}_t\big]\dif t \\
    &\quad+\tilde\mE\big[\p_\nu b(Y^{y,\nu}_t,\cL_{Y^\eta_t})(\tilde Y^{\tilde \eta}_t)\cdot \tilde Z^{\tilde \eta}_t(\tilde y)\big]\dif t+\p_y \sigma(Y^{y,\nu}_t,\cL_{Y^\eta_t})\cdot\p_\nu Y^{y,\nu}_t(\tilde y) \dif W_t\\
    &\quad+\tilde\mE\big[\p_\nu \sigma(Y^{y,\nu}_t,\cL_{Y^\eta_t})(\tilde Y^{\tilde y,\nu}_t)\cdot\p_y \tilde Y^{\tilde y,\nu}_t\big]\dif W_t \\
    &\quad+\tilde\mE\big[\p_\nu \sigma(Y^{y,\nu}_t,\cL_{Y^\eta_t})(\tilde Y^{\tilde \eta}_t)\cdot \tilde Z^{\tilde \eta}_t(\tilde y)\big]\dif W_t,
	\end{align*}
	and
	\begin{align*}
    \dif Z^{\eta}_t(\tilde y)&=\p_y b(Y^{\eta}_t,\cL_{Y^\eta_t})\cdot Z^{\eta}_t(\tilde y)\dif t+\tilde\mE\big[\p_\nu b(Y^{\eta}_t,\cL_{Y^\eta_t})(\tilde Y^{\tilde y,\nu}_t)\cdot\p_y \tilde Y^{\tilde y,\nu}_t\big]\dif t \\
    &\quad+\tilde\mE\big[\p_\nu b(Y^{\eta}_t,\cL_{Y^\eta_t})(\tilde Y^{\tilde \eta}_t)\cdot \tilde Z^{\tilde \eta}_t(\tilde y)\big]\dif t+\p_y \sigma(Y^{\eta}_t,\cL_{Y^\eta_t})\cdot Z^{\eta}_t(\tilde y) \dif W_t\\
    &\quad+\tilde\mE\big[\p_\nu \sigma(Y^{\eta}_t,\cL_{Y^\eta_t})(\tilde Y^{\tilde y,\nu}_t)\cdot\p_y \tilde Y^{\tilde y,\nu}_t\big]\dif W_t \\
    &\quad+\tilde\mE\big[\p_\nu \sigma(Y^{\eta}_t,\cL_{Y^\eta_t})(\tilde Y^{\tilde \eta}_t)\cdot \tilde Z^{\tilde \eta}_t(\tilde y)\big]\dif W_t.
	\end{align*}
    By It\^{o} formula, we have
    \begin{align*}
	\dif \mE\|Z^{\eta}_t(\tilde y)\|^p
    &\leq\frac{p}{2}\mE\big[\|Z^{\eta}_t(\tilde y)\|^{p-2}\cdot2\langle Z^{\eta}_t(\tilde y),\p_y b(Y^{\eta}_t,\cL_{Y^{\eta}_t})\cdot Z^{\eta}_t(\tilde y) \\
    &\quad+\tilde \mE[\p_\nu b(Y^{\eta}_t,\cL_{Y^{\eta}_t})(\tilde Y^{\tilde \eta}_t)\cdot\tilde Z^{\tilde \eta}_t(\tilde y)]
    +\tilde \mE[\p_\nu b(Y^{\eta}_t,\cL_{Y^{\eta}_t})(\tilde Y^{\tilde y,\nu}_t)\cdot\p_y\tilde Y^{\tilde y,\nu}_t]\rangle\big]\dif t\\
    &\quad+\frac{p(p-1)}{2}\mE\big[\|Z^{\eta}_t(\tilde y)\|^{p-2}\cdot\|\p_y \sigma(Y^{\eta}_t,\cL_{Y^{\eta}_t})\cdot Z^{\eta}_t(\tilde y) \\
    &\quad+\tilde \mE[\p_\nu \sigma(Y^{\eta}_t,\cL_{Y^{\eta}_t})(\tilde Y^{\tilde \eta}_t)\cdot\tilde Z^{\tilde \eta}_t(\tilde y)]
    +\tilde \mE[\p_\nu \sigma(Y^{\eta}_t,\cL_{Y^{\eta}_t})(\tilde Y^{\tilde y,\nu}_t)\cdot\p_y\tilde Y^{\tilde y,\nu}_t]\|^2\big]\dif t.
	\end{align*}
    By the assumption {\bf ($\hat H^{\sigma,b}$)}, we have for any $h\in\mR^{d_2}$ and $H\in L^2(\Omega)$,
	\begin{align*}
    &2\langle h,\p_yb(y,\nu)\cdot h+\mE[\p_\nu b(y,\nu)(\eta)\cdot H]\rangle+(p-1)\|\p_y\sigma(y,\nu)\cdot h+\mE[\p_\nu\sigma(y,\nu)(\eta)\cdot H]\|^2\\
	&\leq c_1\mE|H|^2-c_2|h|^2,
	\end{align*}
	with $\cL_\eta=\nu$. Consequently, we arrive at
    \begin{align*}
     \dif \mE\|Z^{\eta}_t(\tilde y)\|^p&\leq-\frac{p}{2}(c_2-c_1-\gamma)\mE\|Z^{\eta}_t(\tilde y)\|^p\dif t+C_0\mE\|\p_yY^{\tilde y,\nu}_t\|^p\dif t\\
     &\leq-\frac{p}{2}(c_2-c_1-\gamma)\mE\|Z^{\eta}_t(\tilde y)\|^p\dif t+C_0e^{-\frac{p}{2}c_2t}\dif t,
    \end{align*}
    which together with the comparison theorem yields
    \begin{align*}
     \mE\|Z^{\eta}_t(\tilde y)\|^p\leq C_0e^{-\frac{p}{2}(c_2-c_1-\gamma)t}.
    \end{align*}
    In the same way we deduce
    \begin{align*}
	&\dif \mE\|\p_\nu Y^{y,\nu}_t(\tilde y)\|^p\\
    &\leq\frac{p}{2}\mE\big[\|\p_\nu Y^{y,\nu}_t(\tilde y)\|^{p-2}\cdot2\langle \p_\nu Y^{y,\nu}_t(\tilde y),\p_y b(Y^{y,\nu}_t,\cL_{Y^{\eta}_t})\cdot\p_\nu Y^{y,\nu}_t(\tilde y) \\
    &\quad+\tilde \mE[\p_\nu b(Y^{y,\nu}_t,\cL_{Y^{\eta}_t})(\tilde Y^{\tilde \eta}_t)\cdot\tilde Z^{\tilde \eta}_t(\tilde y)]
    +\tilde \mE[\p_\nu b(Y^{y,\nu}_t,\cL_{Y^{\eta}_t})(\tilde Y^{\tilde y,\nu}_t)\cdot\p_y\tilde Y^{\tilde y,\nu}_t(\tilde y)]\rangle\big]\dif t\\
    &\quad+\frac{p(p-1)}{2}\mE\big[\|\p_\nu Y^{y,\nu}_t(\tilde y)\|^{p-2}\cdot\|\p_y \sigma(Y^{y,\nu}_t,\cL_{Y^{\eta}_t})\cdot\p_\nu Y^{y,\nu}_t(\tilde y) \\
    &\quad+\tilde \mE[\p_\nu \sigma(Y^{y,\nu}_t,\cL_{Y^{\eta}_t})(\tilde Y^{\tilde \eta}_t)\cdot\tilde Z^{\tilde \eta}_t(\tilde y)]
    +\tilde \mE[\p_\nu \sigma(Y^{y,\nu}_t,\cL_{Y^{\eta}_t})(\tilde Y^{\tilde y,\nu}_t)\cdot\p_y\tilde Y^{\tilde y,\nu}_t(\tilde y)]\|^2\big]\dif t\\
    &\leq-\frac{p}{2}(c_2-\gamma)\mE\|\p_\nu Y^{y,\nu}_t(\tilde y)\|^p\dif t+C_0\mE\|Z^{\eta}_t(\tilde y)\|^p\dif t+C_0\mE\|\p_yY^{\tilde y,\nu}_t\|^p\dif t\\
     &\leq-\frac{p}{2}(c_2-\gamma)\mE\|\p_\nu Y^{y,\nu}_t(\tilde y)\|^p\dif t+C_0e^{-\frac{p}{2}(c_2-c_1-\gamma)t}\dif t,
	\end{align*}
    and thus
    \begin{align*}
     \mE\|\p_\nu Y^{y,\nu}_t(\tilde y)\|^p\leq C_0e^{-\frac{p}{2}(c_2-c_1-\gamma)t}.
    \end{align*}
    Similarly, we have
    \begin{align*}
	&\dif \mE\|\p_{\tilde y}Z^{\eta}_t(\tilde y)\|^p\\
    &\leq-\frac{p}{2}(c_2-c_1-\gamma)\mE\|\p_{\tilde y}Z^{\eta}_t(\tilde y)\|^p\dif t+C_0\mE\|\p_yY^{\tilde y,\nu}_t\|^{2p}\dif t+C_0\mE\|\p^2_yY^{\tilde y,\nu}_t\|^p\dif t\\
     &\leq-\frac{p}{2}(c_2-c_1-\gamma)\mE\|\p_{\tilde y}Z^{\eta}_t(\tilde y)\|^p\dif t+C_0e^{-\frac{p}{2}(c_2-\gamma)t}\dif t,
	\end{align*}
    which in turn yields
    \begin{align*}
     \mE\|\p_{\tilde y}Z^{\eta}_t(\tilde y)\|^p\leq C_0e^{-\frac{p}{2}(c_2-c_1-\gamma)t}.
    \end{align*}
    This further implies
    \begin{align*}
    \dif \mE\|\p_{\tilde y}\p_\nu Y^{y,\nu}_t(\tilde y)\|^p
    &\leq-\frac{p}{2}(c_2-\gamma)\mE\|\p_{\tilde y}\p_\nu Y^{y,\nu}_t(\tilde y)\|^p\dif t+C_0\mE\|\p_{\tilde y}Z^{\eta}_t(\tilde y)\|^p\dif t\\
    &\quad+C_0\mE\|\p_yY^{\tilde y,\nu}_t\|^{2p}\dif t+C_0\mE\|\p^2_yY^{\tilde y,\nu}_t\|^p\dif t\\
    &\leq-\frac{p}{2}(c_2-\gamma)\mE\|\p_{\tilde y}\p_\nu Y^{y,\nu}_t(\tilde y)\|^p\dif t+C_0e^{-\frac{p}{2}(c_2-c_1-\gamma)t}\dif t,
	\end{align*}
	and thus the desired result is obtained. In the same way, we can prove the estimates of $\mE\|\p_y\p_\nu Y^{y,\nu}_t(\tilde y)\|^p$ and
$\mE\|\p^2_\nu Y^{y,\nu}_t(\tilde y,\breve{y})\|^p$, we omit the details here.
\end{proof}


\bigskip
\begin{thebibliography}{2}

\bibitem{APV} Abdulle A., Pavliotis G. A. and Vaes U.: Spectral methods for multiscale stochastic differential equations. {\it SIAM/ASA J. Uncertain. Quantif.}, {\bf 5} (2017), 720--761.

\bibitem{BR}V. Barbu and M. R\"ockner: From non-linear Fokker-Planck equations to solutions of distribution dependent SDE. {\it Ann. Probab.}, {\bf 48} (2020), 1902--1920.
	
	\bibitem{BRR}V. Barbu, M. R\"ockner and F. Russo: Probabilistic representation for solutions
	of an irregular porous media type equation: the degenerate case. {\it Probab. Theory Rel. Fields}, {\bf 15} (2011), 1--43.
	
	

	
	
		
		\bibitem{BS2}Z. W. Bezemek and K. Spiliopoulos:  Large deviations for interacting multiscale particle systems. arXiv: 2011.03032.
	
	\bibitem{BS}
	Z. W. Bezemek and K. Spiliopoulos: Rate of homogenization for fully-coupled McKean-Vlasov SDEs. arXiv:2202.07753v1.
	

	
	
	\bibitem{BS3}Z. W. Bezemek and K. Spiliopoulos: Moderate deviations for fully coupled multiscale weakly interacting particle systems. arXiv: 2202.08403.
	
\bibitem{BSV}	V. I. Bogachev, S. V. Shaposhnikov and A. Yu. Veretennikov: Differentiability of solutions of stationary Fokker-Planck-Kolmogorov equations with respect to a parameter. {\it Discrete   Cont. Dynam. Syst.},  {\bf  36} (2016),  3519--3543.
	
	\bibitem{BLPR} R. Buckdahn, J. Li, S. Peng and C. Rainer: Mean-field stochastic differential equations and associated PDEs. {\it Ann. Probab.}, {\bf 45} (2017), 824--787.

\bibitem{C} P. Cardaliaguet P.:  Notes on mean field games. https:/\!/www.ceremade.dauphine.fr/cardaliaguet/MF-G20130420.pdf, 2013.

\bibitem{CGT}E. A. Carlen, E. Gabetta and G. Toscani: Propagation of smoothness and the rate of exponential convergence to equilibrium for a spatially homogeneous Maxwellian gas. {\it Commun. Math. Phys.}, {\bf 305} (1999), 521--546.
	
	\bibitem{CD}
	R. Carmona and F. Delarue: Probabilistic Theory of Mean Field Games with Applications I: Mean Field FBSDEs, Control, and Games, Probability Theory and Stochastic Modelling, Springer, 2018.
	
	\bibitem{CC2}J. A. Carrillo and  Y. P. Choi: Quantitative error estimates for the large friction limit of Vlasov equation with nonlocal forces. {\it Ann. Ins. Henri Poincar\'e, Anal. non lin.}, {\bf37} (2020), 925--954.

\bibitem{CC}J. A. Carrillo and  Y. P. Choi: Mean-field limits: from particle descriptions to macroscopic equations. {\it Arch. Rat. Mech. Anal.}, {\bf 241} (2021), 1529--1573.
	
	\bibitem{CGPS}J. A. Carrillo, R. S. Gvalani, G. A. Pavliotis and A. Schlichting: Long-time behaviour and phase transitions for the Mckean-Vlasov equation on the torus. {\it Arch. Rational Mech. Anal.}, {\bf 235} (2020), 635--690.


\bibitem{CCD}J. F. Chassagneux, D. Crisan and F. Delarue. A probabilistic approach to classical solutions of the master equation for large population equilibria. to appear in {\it Memoirs of the AMS},	arXiv:1411.3009.
	
	\bibitem{CF} P.-E. Chaudru de Raynal and N. Frikha: Well-Posedness for some non-linear SDEs and related PDE on the Wassertein space. {\it J. Math. Pures Appl.}, {\bf 159} (2022), 1--167.

\bibitem{CCKW} X. Chen, Z.-Q. Chen, T. Kumagai and J. Wang: Periodic homogenization of non-symmetric L\'evy-type processes. {\it Ann. Probab.}, {\bf 49} (2021), 2874--2921.
	
	

	
	

	
	
	\bibitem{CM}D. Crisan and E. McMurray: Smoothing properties of McKean-Vlasov SDEs. {\it Probab. Theory and Related Fields}, {\bf 171} (2018),  97--148.
	
		\bibitem{D}D. A. Dawson: Critical dynamics and fluctuations for a mean-field model of cooperative behavior. {\it J. Stat. Phys.}, {\bf 31} (1983), 29--85.
	
	\bibitem{GGP}M. G. Delgadino, R. S. Gvalani and G. A. Pavliotis: On the diffusive-mean field limit for weakly interacting diffusions exhibiting phase transitions.  {\it Arch. Rational Mech. Anal.}, {\bf 241} (2021), 91--148.

\bibitem{DGPS}M. G. Delgadino, R. S. Gvalani, G. A. Pavliotis and S. A. Smith: Phase transitions, logarithmic Sobolev inequalities, and uniform-in-time propagation of chaos for weakly interacting diffusions. arXiv:2112.06304v1.

\bibitem{DGO}    M. Duerinckx, A. Gloria and F. Otto: The structure of fluctuations in stochastic homogenization. {\it Commun. Math. Phy.}, {\bf 377} (2020), 259--306.
	
	
	\bibitem{DP}    A. B. Duncan and G. A. Pavliotis: Brownian motion in an $N$-scale periodic potential. arXiv:1605.05854.
	
	\bibitem{DEGS}A. Durmus, A. Eberle, A Guillin and K Schuh:	Sticky nonlinear SDEs and convergence of McKean-Vlasov equations without confinement. arXiv:2201.07652.

\bibitem{EGZ}A. Eberle, A. Guillin and R. Zimmer: Quantitative Harris type theorems for diffusions and McKean-Vlasov processes. {\it Trans. Amer. Math. Soc.}, {\bf 371} (2019), 7135--7173.


\bibitem{FG}N. Fournier and A. Guillin: From a Kac-like particle system to the Landau equation for hard potentials and Maxwell molecules. {\it Ann. Sci. Éc. Norm. Supér.},  {\bf 50} (2017), 157--199.
	
	\bibitem{GP}S. N. Gomes and G. A. Pavliotis: Mean field limits for interacting diffusions in a two-scale
	potential. {\it J. Non-linear Sci.}, {\bf 28} (2018), 905-941.

\bibitem{GS} R. S. Gvalania and A. Schlichting: Barriers of the McKean-Vlasov energy via a mountain pass theorem in the space of probability measures. {\it J. Func. Anal.}, {\bf 279} (2020), 108720.


\bibitem{HP}M. Hairer and E. Pardoux: Fluctuations around a homogenised semilinear random
PDE. {\it Arch. Ration. Mech. Anal.}, {\bf 239} (2021), 151--217.
	
	

	

	
	\bibitem{HLL2}
	W. Hong, S. Li and W. Liu: Strong convergence rates in averaging principle for slow-fast McKean-Vlasov SPDEs. arXiv:2107.14401v1.
	
	\bibitem{HMP}	I. Honor\'e, S. Menozzi and G. Pag\`es: Non-asymptotic Gaussian estimates for the recursive approximation of the invariant measure of a diffusion. {\it Ann. Inst Henri Poincare-Pr.}, {\bf 56} (2020),  1559--1605.
	
	
	\bibitem{HW}
	X. Huang and F.-Y. Wang: Distribution dependent SDEs with singular coefficients. {\it Stochastic Process. Appl.}, {\bf129} (2019),  4747--4770.
	
	
	
	

\bibitem{J}P.-E. Jabin: Macroscopic limit of Vlasov type equations with friction. {\it Ann. Ins. Henri Poincar\'e, Anal. non lin.}, {\bf17} (2000), 651--672.

	
\bibitem{K}	M. Kac: Foundations of kinetic theory. In Proceedings of the Third Berkeley Symposium on Mathematical Statistics and Probability: Contributions to Astronomy and Physics, pages 171–197, Berkeley, Calif., 1956. University of California Press.

\bibitem{K2}V. N. Kolokoltsov: Nonlinear Markov processes and kinetic equations. Cambridge University Press, 2010.


	
	
	
\bibitem{L} P. Lions: Mean-field games and applications. Lectures at the College de France, 2007.

\bibitem{MST}J. C. Mattingly, A. M. Stuart and M. V. Tretyakov: Convergence of numerical time-averaging and stationary measures via Poisson equations. {\it SIAM J. Numer. Anal.},  {\bf 48} (2010), 552--577.

	
	
	\bibitem{M}
	H. P. McKean, Jr.: A class of Markov processes associated with nonlinear parabolic equations. {\it Proc. Nat. Acad. Sci. USA}, {\bf56 } (1966), 1907--1911.
	
\bibitem{MV}	Y. S. Mishura and A. Yu. Veretennikov: Existence and uniqueness theorems for solutions of McKean-Vlasov stochastic equations. {\it Theory Probab. Math. Stat.}, {\bf 103} (2021), 59--101.

\bibitem{PP}Pag\`es G. and Panloup F.: Ergodic approximation of the distribution of a stationary diffusion: rate of convergence. {\it Ann. Appl. Probab.}, {\bf 22} (2012), 1059--1100.
	
	\bibitem{P-V} E. Pardoux and A. Yu. Veretennikov: On the Poisson equation and diffusion approximation. I. {\it Ann. Prob.}, {\bf 29} (2001), 1061--1085.
	
	\bibitem{P-V2} E. Pardoux and A. Yu. Veretennikov: On the Poisson equation and diffusion approximation 2. {\it Ann. Prob.}, {\bf 31} (2003), 1166--1192.
	
	\bibitem{P-V3}E. Pardoux and A. Yu. Veretennikov: On the Poisson equation and diffusion approximation 3. {\it Ann. Prob.}, {\bf 33} (2005), 1111--1133.
	
	\bibitem{PS}G. A. Pavliotis and A. M. Stuart: Multiscale Methods. vol. 53. Texts in Applied Mathematics Springer, New York, 2008.
	

	
	\bibitem{RSX}
	M. R\"{o}ckner, X. Sun and Y. C. Xie: Strong convergence order for slow-fast McKean-Vlasov stochastic differential equations. {\it Ann. Inst. H. Poincar\'{e} Probab. Statist.}, {\bf57} (2021), 547--576.
	
	
	
	\bibitem{RX1} M. R\"ockner and L. Xie: Diffusion approximation for fully coupled stochastic differential equations. {\it Ann. Prob.}, {\bf 49} (2021), 1205--1236.
	

	
	
	\bibitem{W1}F.-Y. Wang: Distribution dependent SDEs for Landau type equations. {\it Stoch. Proc. Appl.}, {\bf 128} (2018), 595--621.
	
	
\end{thebibliography}
\end{document}